\input jytex.tex   
\typesize=10pt \magnification=1200 \baselineskip17truept
\footnotenumstyle{arabic} \hsize=6truein\vsize=8.5truein
\sectionnumstyle{blank}
\chapternumstyle{blank}
\chapternum=1
\sectionnum=1
\pagenum=0

\def\begintitle{\pagenumstyle{blank}\parindent=0pt
\begin{narrow}[0.4in]}
\def\endtitle{\end{narrow}\newpage\pagenumstyle{arabic}}


\def\beginexercise{\vskip 20truept\parindent=0pt\begin{narrow}[10
truept]}
\def\endexercise{\vskip 10truept\end{narrow}}


\def\eql#1{\eqno\eqnlabel{#1}}
\def\ref{\reference}
\def\peq{\puteqn}
\def\pref{\putref}

\def\mgn{\marginnote}
\def\bex{\begin{exercise}}
\def\eex{\end{exercise}}


\font\open=msbm10 


\def\StretchRtArr#1{{\count255=0\loop\relbar\joinrel\advance\count255 by1
\ifnum\count255<#1\repeat\rightarrow}}
\def\StretchLtArr#1{\,{\leftarrow\!\!\count255=0\loop\relbar
\joinrel\advance\count255 by1\ifnum\count255<#1\repeat}}

\def\StretchLRtArr#1{\,{\leftarrow\!\!\count255=0\loop\relbar\joinrel\advance
\count255 by1\ifnum\count255<#1\repeat\rightarrow\,\,}}

\def\mbox#1{{\leavevmode\hbox{#1}}}

\def\hspace#1{{\phantom{\mbox#1}}}
\def\oR{\mbox{\open\char82}}

\def\oZ{\mbox{\open\char90}}

\def\oQ{\mbox{\open\char81}}

\def\al{\alpha}
\def\bmu{{\bmit\mu}} 
\def\be{\beta}
\def\ga{\gamma}

\def\Ga{\Gamma}

\def\om{\omega}

\def\th{\theta}

\def\ze{\zeta}

\def\De{\Delta}

\def\caC{{\cal C}}

\def\zf{$\zeta$--function}
\def\zfs{$\zeta$--functions}


\def\frac#1/#2{\leavevmode\kern.1em
\raise.5ex\hbox{\the\scriptfont0 #1}\kern-.1em/\kern-.15em
\lower.25ex\hbox{\the\scriptfont0 #2}}
\def\sfrac#1/#2{\leavevmode\kern.1em
\raise.5ex\hbox{\the\scriptscriptfont0 #1}\kern-.1em/\kern-.15em
\lower.25ex\hbox{\the\scriptscriptfont0 #2}}

\def\gtorder{\mathrel{\raise.3ex\hbox{$>$}\mkern-14mu
             \lower0.6ex\hbox{$\sim$}}}
\def\ltorder{\mathrel{\raise.3ex\hbox{$<$}\mkern-14mu
             \lower0.6ex\hbox{$\sim$}}}

\def\semidirprod{\rlap{\ss C}\raise1pt\hbox{$\mkern.75mu\times$}}
\def\for{\lower6pt\hbox{$\Big|$}}
\def\fish{\kern-.25em{\phantom{abcde}\over \phantom{abcde}}\kern-.25em}


\def\boxit#1{\vbox{\hrule\hbox{\vrule\kern3pt
        \vbox{\kern3pt#1\kern3pt}\kern3pt\vrule}\hrule}}
\def\dalemb#1#2{{\vbox{\hrule height .#2pt
        \hbox{\vrule width.#2pt height#1pt \kern#1pt \vrule
                width.#2pt} \hrule height.#2pt}}}

\def\ol{\overline}
\def\frac#1#2{{{#1}\over{#2}}}

\def\noin{\noindent}


\def\eg{{\it e.g.}}
\def\ie{{\it i.e. }}
\def\cf{{\it cf }}
\def\pa{\partial}


\def\Tr{{\rm Tr\,}}

\def\wt{\widetilde}

\def\3j#1#2#3#4#5#6{\left\lgroup\matrix{#1&#2&#3\cr#4&#5&#6\cr}
\right\rgroup}

\def\man{{\cal M}}

\def\m?{\mgn{?}}

\def\pa{\partial}

\def\beq{\begin{eqnarray}}
\def\eeq{\end{eqnarray}}


\def\aop#1#2#3{{\it Ann. Phys.} {\bf {#1}} ({#2}) #3}

\def\cmp#1#2#3{{\it Comm. Math. Phys.} {\bf {#1}} ({#2}) #3}
\def\cqg#1#2#3{{\it Class. Quant. Grav.} {\bf {#1}} ({#2}) #3}

\def\ijmp#1#2#3{{\it Int. J. Mod. Phys.} {\bf {#1}} ({#2}) #3}

\def\jmp#1#2#3{{\it J. Math. Phys.} {\bf {#1}} ({#2}) #3}
\def\jpa#1#2#3{{\it J. Phys.} {\bf A{#1}} ({#2}) #3}
\def\lnm#1#2#3{{\it Lect. Notes Math.} {\bf {#1}} ({#2}) #3}

\def\np#1#2#3{{\it Nucl. Phys.} {\bf B{#1}} ({#2}) #3}
\def\pl#1#2#3{{\it Phys. Lett.} {\bf {#1}} ({#2}) #3}

\def\prp#1#2#3{{\it Phys. Rep.} {\bf {#1}} ({#2}) #3}
\def\pr#1#2#3{{\it Phys. Rev.} {\bf {#1}} ({#2}) #3}
\def\prA#1#2#3{{\it Phys. Rev.} {\bf A{#1}} ({#2}) #3}

\def\prD#1#2#3{{\it Phys. Rev.} {\bf D{#1}} ({#2}) #3}
\def\prl#1#2#3{{\it Phys. Rev. Lett.} {\bf #1} ({#2}) #3}

\def\rmp#1#2#3{{\it Rev. Mod. Phys.} {\bf {#1}} ({#2}) #3}

\def\zfp#1#2#3{{\it Z. f. Phys.} {\bf {#1}} ({#2}) #3}

\def\cras#1#2#3{{\it Comptes Rend. Acad. Sci. (Paris)} {\bf{#1}} (#2) #3}
\def\prs#1#2#3{{\it Proc. Roy. Soc.} {\bf A{#1}} ({#2}) #3}
\def\pcps#1#2#3{{\it Proc. Camb. Phil. Soc.} {\bf{#1}} ({#2}) #3}
\def\mpcps#1#2#3{{\it Math. Proc. Camb. Phil. Soc.} {\bf{#1}} ({#2}) #3}

\def\amsh#1#2#3{{\it Abh. Math. Sem. Ham.} {\bf {#1}} ({#2}) #3}
\def\am#1#2#3{{\it Acta Mathematica} {\bf {#1}} ({#2}) #3}
\def\aim#1#2#3{{\it Adv. in Math.} {\bf {#1}} ({#2}) #3}
\def\ajm#1#2#3{{\it Am. J. Math.} {\bf {#1}} ({#2}) #3}

\def\aom#1#2#3{{\it Ann. of Math.} {\bf {#1}} ({#2}) #3}
\def\cjm#1#2#3{{\it Can. J. Math.} {\bf {#1}} ({#2}) #3}
\def\bams#1#2#3{{\it Bull.Am.Math.Soc.} {\bf {#1}} ({#2}) #3}

\def\cmh#1#2#3{{\it Comm. Math. Helv.} {\bf {#1}} ({#2}) #3}

\def\dmj#1#2#3{{\it Duke Math. J.} {\bf {#1}} ({#2}) #3}
\def\invm#1#2#3{{\it Invent. Math.} {\bf {#1}} ({#2}) #3}

\def\jdg#1#2#3{{\it J. Diff. Geom.} {\bf {#1}} ({#2}) #3}

\def\joa#1#2#3{{\it J. of Algebra} {\bf {#1}} ({#2}) #3}
\def\jram#1#2#3{{\it J. f. reine u. Angew. Math.} {\bf {#1}} ({#2}) #3}
\def\jims#1#2#3{{\it J. Indian. Math. Soc.} {\bf {#1}} ({#2}) #3}
\def\jlms#1#2#3{{\it J. Lond. Math. Soc.} {\bf {#1}} ({#2}) #3}
\def\jmpa#1#2#3{{\it J. Math. Pures. Appl.} {\bf {#1}} ({#2}) #3}
\def\ma#1#2#3{{\it Math. Ann.} {\bf {#1}} ({#2}) #3}

\def\mz#1#2#3{{\it Math. Zeit.} {\bf {#1}} ({#2}) #3}
\def\ojm#1#2#3{{\it Osaka J.Math.} {\bf {#1}} ({#2}) #3}

\def\pems#1#2#3{{\it Proc. Edin. Math. Soc.} {\bf {#1}} ({#2}) #3}

\def\plb#1#2#3{{\it Phys. Letts.} {\bf {B#1}} ({#2}) #3}
\def\pla#1#2#3{{\it Phys. Letts.} {\bf {A#1}} ({#2}) #3}
\def\plms#1#2#3{{\it Proc. Lond. Math. Soc.} {\bf {#1}} ({#2}) #3}
\def\pgma#1#2#3{{\it Proc. Glasgow Math. Ass.} {\bf {#1}} ({#2}) #3}
\def\qjm#1#2#3{{\it Quart. J. Math.} {\bf {#1}} ({#2}) #3}
\def\qjpam#1#2#3{{\it Quart. J. Pure and Appl. Math.} {\bf {#1}} ({#2}) #3}

\def\rmjm#1#2#3{{\it Rocky Mountain J. Math.} {\bf {#1}} ({#2}) #3}

\def\tams#1#2#3{{\it Trans.Am.Math.Soc.} {\bf {#1}} ({#2}) #3}

\begin{title}
\vglue 0.5truein
\vskip15truept
\centertext {\Bigfonts \bf Analytic torsion on spherical factors}
\vskip7truept

 \centertext {\Bigfonts \bf and tessellations}
\vskip10truept \centertext{\Bigfonts \bf }
 \vskip 20truept
\centertext{J.S.Dowker\footnote{dowker@man.ac.uk}} \vskip 7truept
\centertext{\it Theory Group,} \centertext{\it School of Physics and
Astronomy,} \centertext{\it The University of Manchester,}
\centertext{\it Manchester, England} \vskip 7truept \centertext{and}
\vskip 7truept

 \centertext{Peter
Chang} \vskip 7truept \centertext{\it Diamond Light Source Ltd.,}
\centertext{\it Harwell, Oxon., UK.}

\vskip40truept
\begin{narrow}
The analytic torsion is computed on fixed--point free and non
fixed--point free factors (tessellations) of the three--sphere. We repeat
the standard computation on spherical space forms (Clifford--Klein
spaces) by an improved technique. The transformation to a simpler form of
the spectral expression of the torsion on spherical factors effected by
Ray is shown to be more general than his derivation implies. It
effectively allows the eigenvalues to be considered as squares of
integers, and applies also to trivial twistings. The analytic torsions
compute to algebraic numbers, as expected. In the case of icosahedral
space, the quaternion twisting gives a torsion proportional to the
fundamental unit of $\oQ(\sqrt5)$. As well as a direct calculation, the
torsions are obtained from the lens space values by a character inducing
procedure. On tessellations, terms occur due to edge conical
singularities.
\end{narrow}
\vskip 5truept
\vskip 60truept
\vfil
\end{title}
\pagenum=0
\newpage

\section{\bf 1. Introduction.}

The Reidemeister torsion was introduced originally as a secondary
combinatorial topological invariant and used by Franz to distinguish
between non-- homeomorphic lens spaces that were cohomologically and
homotopically identical. An analytical analogue was later defined by Ray
and Singer in the Riemannian de Rham setting and shown by Cheeger and
M\"uller to be equal to the Reidemeister torsion. There is a useful
summary in the nice little book by Rosenberg, [\pref{Rosenberg}]. The
evaluation for lens spaces is classic, [\pref{Ray}].

The torsion has become of some physical significance following the work
of Schwarz on topological field theory, and a physicist's calculation of
the torsion on lens spaces can be found in Nash and O'Connor,
[\pref{NandO}].

It was initially defined on manifolds without boundary, although Cheeger
extended it to those with a boundary. Later, Lott and Rothenberg, L\"uck
and Vishik showed that, in this case, when the metric is product at the
boundary, there is a difference between the Reidemeister and Ray--Singer
quantities which involved the boundary Euler number.

At around the same time, the second named author calculated the torsion
on the tessellated three--sphere using the information gathered during  a
computation of the Casimir energy on orbifold factors of spheres and the
results were given in his thesis of 1993, [\pref{Chang}]. The present
paper details these results and adds a few up--to--date comments in the
light of continuing interest in analytic torsion, \eg\ Vertman,
[\pref{Vertman}], and De Melo {\it et al}, [\pref{MHS}]. We particularly
wish to draw attention to the cancellation (section 4) which improves a
result of Ray, [\pref{Ray}].
\section{\bf 2. The tessellation.}

The tessellation is classic. In related calculations it has been
described in our earlier works, [\pref{ChandD}], [\pref{Dow20}], where
many original references are given. In the orbifold factoring, S$^3/\Ga$,
non--trivial elements of the full reflective, polytope symmetry group,
$\Ga$, leave hyperplanes, planes, lines and points invariant (in
$\oR^4$), and the fundamental domain (a spherical tetrahedron), $\man$,
on S$^3$, has a boundary, as well as edges and vertices. Restricting to
the rotational subgroup, $\Ga^+$, doubles the fundamental domain and
removes the odd co--dimension singular domains, in this case this means
the two--dimensional boundary, $\pa\man$, and the vertices. The remaining
singular edges form the axes of periodic dihedral wedges (cones).
\section{\bf 3. Analytic torsion.}
We need a few, very standard formulae and preliminary results. The
definition of analytic torsion, $T$, by Ray and Singer reads, for a
manifold, $\man$, of dimension $d$,
  $$\eqalign{
   \log T(\man,\rho)&={1\over2}\sum_{p=0}^d p\,(-1)^p\,\ze'_{\De_p,\rho}(0)\cr
   }
   \eql{at1}
  $$
where $\ze_{\De_p,\rho}(s)$ is the \zf\ for the de Rham Laplacian,
twisted by a representation, $\rho$, of $\Ga$, \ie
  $$
  (\ga\phi)(x)=\phi(\ga x)=\rho(\ga)\phi(x)\,,
  \quad\forall\ga\in\Ga\,,\quad x\in\man\,.
  \eql{twist}
  $$
The form $\phi$ takes values in the flat vector bundle associated with
the representation $\rho$ which is taken to be either orthogonal, $\Ga\to
$O$(N)$, or unitary, $\Ga\to$ U$(N)$. In the former case the form is
assumed real and in the latter, complex. For complex forms, one must be
careful to include a factor of two when counting dimensions.

If $x$ is a fixed point, $\ga x=x$, and, if the rep is non--trivial, this
implies that the form, $\phi$, has to vanish at the fixed points \ie, in
the present situation, on the axis of the cone. In this case, normal
Hodge duality applies. There is no need for absolute or relative
conditions in the purely rotational case which concerns us here.

If the rep, $\rho$, is trivial, then one has to decide on the conditions
that hold at the edge singularity. We select, by default, those that
yield the Friedrichs extension.

An alternative, perhaps more fundamental, form for the torsion is given
in terms of the {\it coexact} \zfs, $\ze_{p,\,\rho}$,
 $$
 \log T(\man,\rho)={1\over2}\sum_{p=0}^{d-1} (-1)^p\,\ze'_{p,\rho}(0)
 $$
which transcribes to (\peq{at1}) using t\'elescopage via the well--known
relation,
  $$
  \ze_{\De_p,\,\rho}(s)=\ze_{p,\,\rho}(s)+\ze_{p-1,\,\rho}(s)\,,
  $$
between the total \zf\ and the coexact \zfs.

At this point we restrict to odd dimensional manifolds, and set $d=2M+1$.
Then (\peq{at1}), for example, is easily rewritten as (see also
[\pref{Nando2}]),
  $$\eqalign{
   \log T(\man,\rho)&=-\sum_{p=0}^{M-1}(-1)^p\,\ze'_{p,\,\rho}(0)
   +{1\over2}(-1)^{M+1}\ze'_{M,\,\rho}(0)\cr
   }
   \eql{at2}
  $$
using Hodge duality and t\'elescopage.

In particular, in the three--dimensional case,
  $$\eqalign{
   \log T(\man,\rho)&={1\over2}\,\ze'_{1,\,\rho}(0)-\ze'_{0,\,\rho}(0)\cr
   }
   \eql{at3}
  $$
This is the case computed in [\pref{Chang}] and now exposed.

To begin with, we define an intermediate quantity,
  $$
  \tau_\rho(s)={1\over2}\,\ze_{1,\,\rho}(s)-\ze_{0,\,\rho}(s)\,.
  \eql{tau}
  $$

In three dimensions, one does not need the full apparatus of $p$--forms.
A coexact one--form is a conformally coupled divergenceless (transverse)
vector and the zero--form is a minimal scalar. The eigenvalues on the
(unit) three-sphere are well--known and the \zfs\ on the orbifold factors
are, (\cf\ [\pref{DandJ}]),
  $$\eqalign{
  \ze_{0,\rho}(s)&=\sum_{n=1}^\infty {d_0(n,\rho)\over(n(n+2))^s}\cr
  \ze_{1,\rho}(s)&=\sum_{n=1}^\infty {d_1(n,\rho)\over(n+1)^{2s}}\,.
  }
  \eql{zetas}
  $$
We remark that any zero mode ($n=0$) that might exist for the scalar, has
been omitted from the sum. This would be the case on the complete sphere,
S$^d$, ($d>1$). For non--trivial twistings, there are no zero modes,
corresponding to a trivial cohomology for (\peq{twist}).

The problem is to find the degeneracies, $d_p$, and then continue the
expression (\peq{tau}). This can be done for the \zfs, (\peq{zetas}),
separately, as in [\pref{NandO}], or for the combination, as in
[\pref{Ray}]. We choose the latter.

In the case of the full sphere, the degeneracies, $d_p(n,{\bf 1})$ are
$(2)N\,(n+1)^2$ and $(2)N\,2n(n+2)$, for $p=0$ and $1$, respectively.
(The factor of 2 comes in for a U($N$) bundle.)

On the factored three--sphere, the necessary formulae have been published
previously in similar contexts, [\pref{DandB,Dow13,Dow40}], so we will
not give the full derivations here. The approach is the very standard one
of group averaging (or projection, or symmetry adaptation) in order to
achieve the twisting, (\peq{twist}). It is convenient (but not necessary)
to use the `left--right' SU(2) actions to represent the action of $\Ga^+$
on S$^3$, which is isomorphic to SU(2).

The situation in [\pref{DandJ}] is very close to the one here. The \zfs\
are diagonal in the (flat) vector bundle indices and, combined with the
explicit discussion in [\pref{Dow40}], this leads to the expression for
the degeneracies,
  $$
  d_0(n-1,\rho)= d(L,L)\,,\quad d_1(n,\rho)=d_\rho(L,L+1)+d_\rho(L+1,L)
  \eql{deg2}
  $$
where $n=2L+1$ with the expected group theory degeneracy,
  $$
  d_\rho(L,J)={1\over|\Ga^+|}\sum_{\ga\in\Ga^+}
  \chi^*_\rho(\ga)\chi^{(L)}(\ga_L)\,\chi^{(J)}(\ga_R)
  \eql{deg1}
  $$
in terms of the character, $\chi_\rho$, of the $\rho$ rep and the SU(2)
character,
  $$
  \chi^{(L)}(\ga)\equiv\chi_{2L+1}(\ga)={\sin
  (2L+1)\th\over\sin\th}\,.
  \eql{char1}
  $$
We note that $\chi_\rho^*({\bf 1})=N$ and, in addition, there is the
factor of two for complex forms.

The angle $\th_L$ is the radial coordinate labelling the SU(2) element
$\ga_L$ as a point on S$^3$.

A tactical decision is whether to leave the group average until last or
to try to effect it earlier. In some cases, \eg\ for cyclic groups, the
latter is possible, and preferable, but here we leave it until last and
so, for convenience, define the summands, which could be termed `partial'
or `off--diagonal' quantities, by the formulae
  $$\eqalign{
   \tau_\rho(s)&={1\over|\Ga^+|}\sum_\ga\chi^*_\rho(\ga)\tau(\ga;s)\cr
   d_p(n,\rho)&={1\over|\Ga^+|}\sum_\ga \chi^*_\rho(\ga)d_p(\ga;n)\,.\cr
   }
   \eql{partial}
  $$

Another decision is whether to treat the two terms in (\peq{tau})
separately, as in [\pref{NandO}], or manipulate them together, as in
[\pref{Ray}]. The former is less economic and does produce interesting
incidental information, however we prefer the second way. To this end,
the vector (spin--one) {\it partial} degeneracy, is rewritten, from
(\peq{deg2}), (\peq{deg1}), (\peq{char1}) and (\peq{partial}), as
  $$
  d_1(\ga;n)=\chi_{n+2}(\th_L)\chi_{n+2}(\th_R)
  +\chi_{n}(\th_L)\chi_{n}(\th_R)-4\cos(n+1)\th_L\,\cos(n+1)\th_R
  $$
to look more like the scalar expression,
  $$
  d_0(\ga;n)=\chi_{n}(\th_L)\chi_{n}(\th_R)\,.
  $$
\section{\bf 4. The cancellation.}

In Ray's calculation, [\pref{Ray}], through the mess of summations and
integrations, a miracle occurs on p.125 in which, for the total quantity,
the effective eigenvalues become those appropriate for a conformally
invariant propagation equation, \ie\ perfect squares. In this section we
give our physicist's version of the implied cancellation, [\pref{Chang}],
which is independent of the free/non--free behaviour of the deck group.
Actually, our result goes further than Ray's as he makes a preliminary
group average for a non--trivial twisting thereby removing a particular
set of terms, which then do not have to be continued. The group average
has still to be performed on the remaining terms.

In contrast, our result is before {\it any} group average and refers to
the complete expression. It can therefore be applied to the case of a
trivial twisting, which provides a useful check, in view of Cheeger's
result, [\pref{Cheeger2}], Theorem 8.35 (see later).

Combining the various \zfs\ one gets
  $$\eqalign{
   \tau(\ga;s)&=\sum_{n=1}^\infty {d_1(\ga;n)\over2(n+1)^{2s}}-{d_0(\ga;n)\over
   \big(n(n+2)\big)^s}\cr
   &=\sum_{n=1}^\infty\bigg[
   {\chi_{n+2}(\th_L)\chi_{n+2}(\th_R)\over2(n+1)^{2s}}-
   {\chi_{n+1}(\th_L)\chi_{n+1}(\th_R)\over\big(n(n+2)\big)^s}\cr
   &\hspace{***}+{\chi_{n}(\th_L)\chi_{n}(\th_R)\over2(n+1)^{2s}}-
   2{\cos(n+1)\th_L\cos(n+1)\th_R\over (n+1)^{2s}}\bigg]\cr
   &={1\over2}\sum_{n=2}^\infty\bigg[
   {1\over(n+1)^s}+{1\over(n-1)^{2s}}-{2\over(n^2-1)^s}\bigg]\chi_{n}(\th_L)\chi_{n}(\th_R)\cr
   &\hspace{******}-2\sum_{n=1}^\infty
   {\cos n\th_L\cos n\th_R\over n^{2s}}+{1\over2^{2s+1}}\,.
   }
   \eql{tau2}
  $$
If we now define two auxiliary functions, $F$ and $\wt\tau$, as the
quantities on the two lines in the last equality, so that,
  $$
  \tau(\ga;s)={1\over2}F(\ga;s)+\wt\tau(\ga;s)\,,
  \eql{tt}
  $$
we proceed to show, quite non--rigorously, that there is no contribution
to the torsion from $F(\ga;s)$.

Now our function, $F$, can be written as,
  $$
  F(\ga;s)=\sum_{n=2}^\infty
  \bigg[{(n-1)^{2s}+(n+1)^{2s}-2(n^2-1)^s\over(n^2-1)^{2s}}\bigg]
  \chi_{n}(\th_L)\chi_{n}(\th_R)\,,
  \eql{eff}
  $$
and it is obvious, by inspection\footnote{ Mathematicians would demand a
higher quality of proof. Some further remarks are given in Appendix 1.},
that it vanishes at $s=0$. Also, taking the derivative with respect to
$s$ gives
  $$\eqalign{
  &{dF(\ga;s)\over ds}=\!\!\sum_{n=2}^\infty\!
  \bigg[{2(n^2-1)^s-(n-1)^{2s}-(n+1)^{2s}\over(n^2-1)^{2s}}\bigg]
  2\ln(n^2-1)\chi_{n}(\th_L)\chi_{n}(\th_R)\cr
  &-\!\!2\!\!\sum_{n=2}^\infty\!
  \bigg[{\ln(n^2\!-\!1)(n^2\!-\!1)^s\!-\!
  \ln(n-1)(n-1)^{2s}\!-\!\ln(n+1)(n+1)^{2s}\over(n^2-1)^{2s}}\bigg]
  \chi_{n}(\th_L)\chi_{n}(\th_R)\big)
  }
  $$
which at $s=0$ yields
  $$\eqalign{
  {dF(\ga;s)\over ds}\bigg|_0&=2\sum_{n=2}^\infty\!
  \big[\ln(n^2-1)-\ln(n-1)-\ln(n+1)\big]
  \chi_{n}(\th_L)\chi_{n}(\th_R)\cr
  &=0\,.
  }
  $$
Hence
  $$
  F(\ga;0)=0\,,\quad{\rm and}\quad F'(\ga;0)=0\,.
  \eql{canc}
  $$
and, therefore, to all intents and purposes, the contribution from
$F(\ga;s)$ can be ignored at the point of interest, $s=0$. Hence,
$\wt\tau$ is related to $\tau$ in (\peq{tt}) by
  $$
  \wt\tau(\ga;0)=\tau(\ga;0)\,,\quad  \wt\tau'(\ga;0)=\tau'(\ga;0)
  $$
and we can work with the simpler, effective tau function, $\wt\tau(s)$,
when calculating the analytic torsion which is, according to (\peq{at3})
and (\peq{tau}),
  $$
  \ln T(\man,\rho)=\wt\tau'_\rho(0)\,.
  \eql{tordef}
  $$
\section{\bf 5. The effective $\tau$ function.}

Putting the group average back, we have the effective function,
  $$
   \wt\tau_{\Ga^+}(s,\rho)={1\over|\Ga^+|}\sum_{\ga\in\,\Ga^+}\chi^*_\rho(\ga)
   \wt\tau(\ga;s)
   \eql{modtau}
  $$
and it is now necessary to address this sum. First some old facts are
needed.

The left and right angles $\th_L$ and $\th_R$ are half the `rotation'
angles associated with the left and right SU(2) actions via the
isomorphism SO(3)$\sim$ SU(2)$/Z_2$. (These rotation angles run from $0$
to $4\pi$.) The $\th_L$ and $\th_R$ are related to angles $\al$ and
$\be$, which are the rotation angles in the two planes of $\oR^4$ under
an SO(4) action, by
  $$
  \th_L={1\over2}(\al-\be)\,,\quad{\rm and}\quad
  \th_R={1\over2}(\al+\be)\,,
  $$
both mod $\pi$.

In the usual embedding, any element of S0(4) is conjugate to a block
diagonal $4\times4$ matrix representation of the form (see \eg\
[\pref{TandS}]),
  $$
  R(\ga)=\left(\matrix { \cos\al&-\sin\al&0&0\cr
          \sin\al&\cos\al&0&0\cr
          0&0&\cos\be&-\sin\be\cr
          0&0&\sin\be&\cos\be\cr}\right)\,.
  \eql{block}
  $$

From (\peq{tau2}) the summand in the effective $\tau$ function,
(\peq{modtau}), is a class function and so the sum can be written as one
over conjugacy classes. The group $\Ga^+$ decomposes into classes,
$\caC_p$, as
   $$
  \Ga^+=\bigoplus_p c_p\,\caC_p
   $$
where $c_p$ is the size of the class labeled by $p$. The total effective
function (\peq{modtau}) then reads
  $$\eqalign{
   \wt\tau_{\Ga^+}(s,\rho)&={1\over|\Ga^+|}\sum_p\chi^*_\rho(p)\,c_p\,
   \wt\tau(\caC_p;s)\cr
   &=\sum_p\chi^*_\rho(p)\,\ol c_p\,
   \wt\tau(\caC_p;s)
   }
   \eql{modtau2}
  $$
with obvious notation and where we have introduced the class `density',
$\ol c_p=c_p/|\Ga^+|$.

In this approach, where the group average is performed by hand, this is
the best that can be done.

From (\peq{tau2}), the partial auxiliary function, $\wt\tau$, is given by
  $$\eqalign{
  \wt\tau(\al,\be;s)&=-\sum_{n=1}^\infty {b(\al,\be,n)\over
  2n^{2s}}+{1\over 2^{2s+1}}\cr
  &=\wt\tau(\caC_p;s)
  }
  \eql{aux}
  $$
where
  $$
  b(\al,\be,n)=4\cos n\th_L\,\cos n\th_R=2(\cos n\al+\cos n\be)
  $$
is recognised as the trace of the SO(4) rep, $R(\ga^n)$, for a rotation
$\ga\sim(\al,\be)$. The same quantity occurs in Ray's treatment,
[\pref{Ray}] p.125. The extension to any odd dimensional sphere is clear
at this point simply by extending the block form, (\peq{block}), granted
the cancellation.

The summation term in (\peq{aux}) is essentially what Ray, [\pref{Ray}],
gets for his $f(s;g)$ on p.125.

One now recognises that the effective $\tau$ function, (\peq{aux}),
involves just a sum of two one--dimensional Epstein \zfs\ defined by,
[\pref{Epstein}],
  $$
   Z\bigg|{g\atop h}\bigg|(s)=\sum_{n=-\infty}^\infty
   |n+g|^{-s}e^{2\pi i nh}\,.
   \eql{epzet}
  $$
If $h=0$, there is the relation with the Hurwitz--Lerch \zf,
  $$
  Z\bigg|{g\atop 0}\bigg|(s)=\ze_R(s,g)+\ze_R(s,1-g)
  \eql{hlzet}
  $$
and, if $g=0$, the $n=0$ term is omitted so that
  $$
  Z\bigg|{0\atop\al/2\pi}\bigg|(s)=2\sum_{n=1}^\infty {\cos n\al\over
  n^{s}}\,,
  \eql{epzet2}
  $$
related to polylogarithms and the Lerch--Lipshitz \zf. This has regularly
occurred in these, and other, situations from the earliest times, \eg\
[\pref{DandB}]. When $\al=0$ it is just twice the Riemann \zf.

We thus have
  $$
  \wt\tau(\al,\be;s)=-{1\over2}\,Z\bigg|{0\atop\al/2\pi}\bigg|(2s)
  -{1\over2}\,Z\bigg|{0\atop\be/2\pi}\bigg|(2s)+{1\over 2^{2s+1}}\,.
  $$

Our primary objective is the torsion and so we now move directly to
evaluate $\wt\tau'(\al,\be;0)$,
  $$
  \wt\tau'(\al,\be;0)=-Z'\bigg|{0\atop\al/2\pi}\bigg|(0)
  -Z'\bigg|{0\atop\be/2\pi}\bigg|(0)-\ln2\,.
  \eql{deriv0}
  $$

It is possible to proceed as in Ray, [\pref{Ray}], but we prefer to
streamline the analysis using some earlier results, as employed in
[\pref{Dow13}] in a computation of lens space determinants.

There are a number of ways of proceeding. Here we first note that as $s$
tends to one,
  $$
   Z\bigg|{h\atop 0}\bigg|(s)\to{2\over s-1}-\psi(h)-\psi(1-h)\,,\quad
   h\ne0\,,
   \eql{lim1}
  $$
Use of the functional relation, for $g=0$,
  $$
   Z\bigg|{g\atop h}\bigg|(2s)=\pi^{2s-1/2}{\Ga(s-1/2)\over\Ga(s)}\,
   e^{-2\pi igh}\,Z\bigg|{h\atop -g}\bigg|(1-2s)\,,
   \eql{funce}
  $$
gives, equivalently, a classic formula,
  $$
  Z'\bigg|{0\atop
  h}\bigg|(0)=-\ga-\ln2\pi-{1\over2}\big(\psi(h)+\psi(1-h)\big)\,,\quad
  h\ne0\,,
  \eql{lim2}
  $$
which is ready to be substituted directly into (\peq{deriv0}).

\section{\bf 6. The partial analytic torsion.}

It is now important to remark that, in the groups we consider, the
angles, $(\th_L,\th_R)$, or $(\al,\be)$, are submultiples of $2\pi$.
Therefore, let, generally,
  $$
  {\al\over2\pi}={a\over q}\,\quad{\rm and}\quad
  {\be\over2\pi}={ b\over q}
  \eql{angles}
  $$
where $a$ and $b$ are integers coprime to $q$. This allows one to use
Gauss' famous formula for $\psi(p/q)$, ($p/q\in\oQ$), or better, a
formula that appears during a proof of this relation, [\pref{Jensen}]
p.146, given below, (see also [\pref{AAR}] p.13), $0<p<q$,
  $$
  \psi\big({p\over q}\big)+\psi\big(1-{p\over q}\big)=-2\ga-2\log q+
  2\sum_{k=1}^{q-1}\cos\big({2\pi p k\over q}\big)\log2\sin{\pi k\over q}\,.
  \eql{gff}
  $$
Substituting into (\peq{deriv0}), there results the equivalent formulae,
  $$\eqalign{
   \wt\tau'(\al,\be;0)
   &=\sum_{k=1}^{q-1}\big(\cos k\al+\cos k\be\big)\log2\sin{\pi k\over q}
   +\ln( 2\pi^2/q^2)\cr
   &=2\sum_{k=1}^{q-1}\cos k\th_L\,\cos k\th_R\,\log2\sin{\pi k\over q}
   +\ln( 2\pi^2/q^2)\cr
   &={1\over2}\sum_{k=1}^{q-1}b(\al,\be,k)\,\log2\sin{\pi k\over q}
   +\ln( 2\pi^2/q^2)\cr
   &\equiv\sum_{k=1}^{q-1}b(\al,\be,k)\,T_q(k)+\ln( 2\pi^2/q^2)
   }
   \eql{torsion3}
  $$
with the angles (\peq{angles}).

Strictly speaking, the special cases when either, or both, of $\al$ and
$\be$ are zero should be treated separately. For example, equation
(\peq{epzet2}) gives the classic value,
  $$
  Z'\bigg|{0\atop 0}\bigg|(0)=2\ze_R'(0)=-\ln2\pi\,,
  $$
but it is easy to show using the ancient product,
  $$
  \prod_{k=1}^{q-1}2\sin({\pi k\over q})=q\,,
  \eql{oldprod}
  $$
that (\peq{torsion3}) covers these cases.

Incidentally, in this connection, equation (\peq{lim1}) can be extended
formally to $h=0$, and (\peq{gff}) to $p=0$ and $q$, by regularising
$\psi(0)$ to $\psi(1), =-\ga$.

For example, for the full sphere, where the group action, and average, is
trivial (all points are fixed),
  $$
  \ln T(S^3,{\bf 1})=\wt\tau'(0,0;0)=\ln2\pi^2=\ln|S^3|\,,
  \eql{fullsp}
  $$
for {\it real} forms, \ie\ $\chi_\rho({\bf 1})=1$. This falls into the
result (\peq{torsion3}) on setting $q=1$, when the sum is non--existent.

Weng and You, [\pref{WandY}], have performed a rather involved
calculation of the analytic torsion on an odd sphere

The old result for the circle is,
  $$
  \ln T(S^1,{\bf 1})=\ln|S^1|\,,
  $$
obtained easily from the Riemann \zf. A simple scaling gives for the
factored circle,
  $$
  \ln T(S^1/\oZ_q,{\bf 1})=\ln|S^1/\oZ_q|\,
  \eql{cone}
  $$

\section{\bf 7. The group average. Clifford--Klein spaces.}

In summary, the (logarithm of the) torsion is given by the twisted group
average (\peq{tordef}) with the effective tau--function (\peq{modtau}) or
(\peq{modtau2}) and the expressions (\peq{torsion3}) for
$\wt\tau'(\caC_p,0)$, the class $\caC_p$ being associated with a pair of
angles, $(\al,\be)$.

As an application, we look at the classic case of fixed--point free
actions \ie Clifford--Klein spaces the basic example of which are lens
spaces. The evaluations are taken a little further than we have seen in
the literature and proceed at a simple level of analysis.

We first rederive Ray's formula for the torsion of the lens space,
$L(q;l_1,l_2)$ $=\,$S$^3/(\oZ_q\times\oZ_q)$, which is defined by the
angles
 $$
  {\al\over2\pi}={ p\,\nu_1\over q}\,,\quad {\be\over2\pi}={p\,\nu_2\over q}\,,
  \eql{angles1}
  $$
where $p,\,=0,\ldots,q-1\,,$ labels $\ga$ (or, equivalently, the class
$\caC_p$). $\nu_1$ and $\nu_2$ are (fixed) integers coprime to $q$, with
$l_1$ and $l_2$ their mod $q$ inverses. With these choices, there are no
fixed points.

By an appropriate selection of a $q$-th root of unity, it would be
possible to set $\nu_1=1$, \ie\ $l_1=1$, without loss of generality. Any
pair, $(\nu_1,\nu_2)$, could be reduced to $(1,\nu)$ by multiplying
through by the mod $q$ inverse of $\nu_1$. The simple, one--sided lens
space, $L(q;1,1)$, corresponds to setting $\nu=1$ so that $\th_L=0$,
$\th_R=2\pi p/q$.

Inserting a U(1) twisting (or, equivalently an SO(2) one), $\chi_r(p)$,
the torsion reads \footnote{ The $\cos k\al$ and $\cos k\be$ terms in
(\peq{torsion3}) can each be replaced by an exponential. Set $k\to q-k$.
We have included an overall factor of two to account for the complex
nature of the forms. Ray does not appear to do this overtly, although he
uses a complex line bundle, as Ray and Singer, [\pref{RandS2}], mainly
do. Ray and Singer in [\pref{RandS}] use orthogonal bundles. For the
single real twisting (when $q$ is even), one could remove the factor of
two. },
  $$\eqalign{
   \ln T\big(L(q;l_1,l_2),r\big)&={2\over q}
   \sum_{j=1}^2\sum_{p=0}^{q-1}\sum_{k=1}^{q-1}
   e^{2\pi i rp/q}\,e^{-2\pi ik\nu_jp/q}\,\ln2\sin{\pi k\over q}\cr
   &={2\over q}\sum_{j=1}^2\sum_{p=0}^{q-1}\sum_{k=1}^{q-1}
   \cos({2\pi rp/q})\,\cos({2\pi k\nu_jp/q})\,\ln2\sin{\pi k\over q}
   }
   \eql{lenst}
  $$
where the integer $r$ determines the bundle twisting, the non--triviality
of which ensures that the constant terms in (\peq{torsion3}) go out.

The sum over $p$ implies that $k\nu_j-r$ is a multiple of $q$ or that
$k\nu_j=r\,$ mod $q$ or that $k=r/\nu_j\,$ mod $q\equiv\,rl_j$.\footnote
{ These manipulations can be found in Epstein, [\pref{Epstein}].}

Therefore we get
  $$
  \ln T\big(L(q;l_1,l_2),r\big)=2\sum_j\ln2\sin{\pi rl_j\over q}
  \eql{tlens}
  $$
which is Ray's value, allowing for the different definition of torsion.
The extension to any odd sphere dimension is obvious.

If the U(1) twisting is trivial (set $r=0$), (\peq{lenst}) yields zero
because $k$ never attains 0 or $q$. The torsion then arises from the
average of the constant term in (\peq{torsion3}) or
  $$\eqalign{
  \ln T(S^3/\oZ_q,{\bf 1})&=\ln (2\pi^2/q)-\ln q\cr
  &=\ln (|S^3/\oZ_q|)-\ln q\,,
  }
  \eql{trivlens}
  $$
where, since the representation is real, the complex doubling has not
been invoked.

The case of a trivial real flat bundle, has been considered by Cheeger,
[\pref{Cheeger2}], Theorem 8.35 which says that the combinatorial
Reidemeister torsion is given by,
  $$
   \ln T_R(\man,{\bf 1})=\sum_{p=0}^d(-1)^{p+1}\big(\ln V_p(\man)+
   \ln O_p\big)
   \eql{flat}
  $$
where $O_p=|{\rm Tor}H^p(\man;\oZ)|=|{\rm Tor}H_{p-1}(\man;\oZ)|$ is the
order of the torsion subgroup of $H^p$ and the $V_p$s  are essentially
volumes associated with the real/integer cohomology and reflect the
existence of zero modes.

In the case that only the top and bottom real cohomology is non--trivial,
one has
   $$
   \ln T_R(\man,{\bf 1})=\ln |\man|-\sum_{p=0}^d(-1)^p\ln O_p
   \eql{flat2}
  $$
after noting that $V_0=1/\sqrt{|\man|}$ and $V_d=\sqrt{|\man|}$. (For
example, the normalised $0$--form zero mode is
$1/\sqrt{|\man|}$.)\footnote{ There appears to be a square root misprint
in Cheeger's text. We refer to the elegant paper of Schwarz and Tyupkin,
[\pref{ScandT}], for a physicist's approach to these questions.}

Applied to a homology lens space, $\man\sim L_q$, this expression gives,
[\pref{Cheeger}],
  $$\eqalign{
  \ln T_R(\man)&=\ln|\man|-{d-1\over2}\ln q\,,\quad d=\dim\man\cr
  }
  \eql{trivtor3}
  $$
The explicit results, (\peq{trivlens}) and (\peq{cone}), of course
confirm this formula.

Rosenberg mentions that the torsion suitably defined  for non--exact
complexes, is (see [\pref{Rosenberg}], p.154),
  $$\eqalign{
  \sum_{p=0}^d(-1)^p\ln|{\rm Tor}H^p(\man;\oZ)|
  =\sum_{p=0}^d(-1)^p\ln|{\rm Tor}H_{p-1}(\man;\oZ)|\,.
  }
  $$
This corresponds to just the first term in Cheeger's Theorem 8.35 and is
a topological invariant, being a special case of Ray and Singer's
extension of the torsion to non--trivial cohomology. See [\pref{RandS2}]
section 3.

The other Clifford--Klein spaces can be computed. The work of Tsuchiya,
[\pref{Tsuchiya}], is concerned with these but does not seem to calculate
any specific values, and uses Ray's formulae (see later). The interesting
and more extensive work of Bauer and Furutani, [\pref{BandF}], must also
be mentioned as it deals with higher dimensions and contains explicit
results.

For simplicity, we consider only one--sided (homogeneous) factors,
S$^3/\Ga'$, where $\Ga'$ is a binary polyhedral group. To be precise, we
set $\th_L=0$ \ie $\al=\be=\th_R\equiv\th_\ga$.

Prism spaces are another infinite set of spaces, for which $\Ga'$ is the
binary dihedral group, $D_q'$, of order $4q$. The generator--relation
structure can be written,
  $$
  A^q=B^2=(AB)^2=Q\,,\quad Q^2=E\,,
  $$
and thus $D_q'$ can be formally written as the direct sum
  $$
  D_q'=Z_{2q}\oplus Z_{2q}B\,,
  $$
where $Z_{2q}$ is generated by $A$. To express the binary doubling, one
has $\oZ_{2q}={\bf Z}_q\oplus Q{\bf Z}_q$ where ${\bf Z}_q$ is the SO(3)
cyclic rotation group and $Q$ is a rotation in $\oR^3$ through 2$\pi$ \ie
$\th_Q=\pi$.

For $A^p$, the angles $\th_\ga$ are,
  $$\eqalign{
  \th_\ga&=\pi p/q\,,\quad p=0,\ldots,2q-1\,.\cr
  }
  \eql{angles2}
  $$
The SU(2) angle, $\th_\ga$, has been left to run from $0$ to $2\pi$,
corresponding to an SO(3) rotation from $0$ to $4\pi$. If one wishes
$\th_\ga$ to be restricted to the range $0$ to $\pi$, as a colatitude on
S$^3$ should be, then it can be arranged that,
  $$\eqalign{
  \th_\ga&=\pi p/q\,,\quad p=0,\ldots,q-1\cr
  &=2\pi-\pi p/q\,,\quad p=q,\ldots,2q-1\,.
  }
  $$

For $\ga=A^pB$, \ie those $2q$ elements containing a (binary) dihedral
rotation, $\th_\ga=\pi/2$ for all $\ga$.

It is straightforward to compute the analytic torsion from
(\peq{torsion3}), using a machine if necessary. The simplest case is the
trivial bundle, and we find
  $$
  \ln T_{D_q'}({\bf 1})=\ln(|S^3|/4q)-2\ln2
  \eql{dihedral}
  $$
which agrees with (\peq{flat2}) in view of the homology, [\pref{TandS}],
  $$\eqalign{
  H_1(S^3/D_q')&=\oZ_4\,,\quad q \,\,{\rm odd} \cr
  &=\oZ_2+\oZ_2\,,\quad q\,\,{\rm even}\,.
  }
  $$

To compute the twisted torsion one needs the irreps of $D_q'$, which are
either one--, or two--dimensional. The former are generated by,
  $$\eqalign{
  \chi(A)=(-1)^a\,,\quad\chi(B)=i^b\,,\quad a=0,1\,,\quad
  b=0,1,2,3\,,
  }
  \eql{1d}
  $$
and the latter by,
  $$
  a(A)=\left(\matrix{e^{i\pi(2a+b)/q}&0\cr
  0&e^{-i\pi(2a+b)/q}}\right)\,,\quad a(B)=\left(\matrix{0&1\cr
  (-1)^b&0}\right)\,,\quad b=0,1\,,
  \eql{2d}
  $$
the conditions being that in (\peq{1d}) for real representations,
$b=0,2$, and then, if $q$ is odd, $a$ cannot equal $1$, while for
imaginary reps, $b=1,3$, then $a=0$ is not possible and, for $a=1$, $q$
must be odd.

In (\peq{2d}), we have
  $$\eqalign{
   {\rm even}\,\,\,q\,\,&\left\{\matrix{\,b=0,\,\quad \,\, a=1, \ldots,q/2-1\cr
   \,b=1,\,\quad \,\, a=0, \ldots,q/2-1\cr}\right.\cr
   {\rm odd}\,\,\,q\,\,&\left\{\matrix{\,b=0,\,\quad \,\, a=1, \ldots,(q-1)/2\cr
   \,b=1,\,\quad \,\, a=0, \ldots,(q-3)/2\cr}\right.\,.\cr
   }
  $$
and for the traces one finds,
  $$
  \Tr (A^p B)=0\,,\quad \Tr A^p=2\cos\bigg({2\pi a p\over q}+{\pi b
  p\over q}\bigg)\,.
  $$

Labelling the one--dimensional irreps by $(a,b)$, calculation gives the
analytic torsion,
  $$\eqalign{
  T_{D_q'}(0,2)&=4/q\,,\quad \forall q\cr
  T_{D_q'}(1,1)=T_{D_q'}(1,3)&=2^2\,,\quad q \,\,{\rm odd}\cr
  T_{D_q'}(1,0)=T_{D_q'}(1,2)&=2\,,\quad q \,\,{\rm even}\,.\cr
  }
  $$
The power of 2 is a complex (or SO(2)) dimension effect.

Labelling the two--dimensional irreps again by $(a,b)$, we find, \eg,
  $$\eqalign{
  &T_{D_2'}(0,1)=2\cr
  &T_{D_3'}(0,1)=1\,,\quad T_{D_3'}(1,0)=3\cr
  &T_{D_4'}(1,0)=2\,,\quad T_{D_4'}(0,1)=2-\sqrt2\,,\quad
  T_{D_4'}(1,1)=2+\sqrt2\,.\cr
  }
  $$

The irreps are usually labelled by letters, for example for $D_4'$ the
two-- dimensional reps are $(1,0)=E$, $(0,1)=E_1'$ and $(1,1)=E_2'$, in
the notation of Landau and Lifshitz, [\pref{LandL}].

No attempt will be made here to give an exhaustive list of all the
numerical possibilities.

The remaining groups are $T'$ (octahedral space), $ O'$ (truncated cube
space) and $Y'$ (dodecahedral, or Poincar\'e, space). For the untwisted
cases we find
    $$\eqalign{
  T_{T'}&={1\over3}\,|S^3|/24\cr
  T_{O'}&={1\over2}\,|S^3|/48\cr
  T_{Y'}&=|S^3|/120\,,\cr
  }
  \eql{CK2}
  $$
in agreement with Cheeger's expression, (\peq{flat}), and the known
(co)homology of the Clifford--Klein manifolds, [\pref{TandS}]. Another
way of obtaining these values is given below and a further one in
Appendix 2.

We compute the twisted values to be,
  $$\eqalign{
  T_{T'}({\bf\dot{1}'})&=(3/2)^2\,,
  \quad T_{T'}({\bf3})=2\cr
  T_{T'}({\bf2_s})&=1/2\,,\quad T_{T'}({\bf\dot{2}'_s})
  =2^2\,\quad\cr
  T_{O'}({\bf1'})&=4/3\,,\quad T_{O'}({\bf2})=3/2\cr
   T_{O'}({\bf3})&=2\,,\quad T_{O'}({\bf3'})=1\,,\quad T_{O'}({\bf4_s})= 2\cr
 T_{O'}({\bf2_s})&=(2-\sqrt2)/2\,,\quad T_{O'}({\bf2'_s})=(2+\sqrt2)/2\cr
  T_{Y'}({\bf 4})&=5/3\,,\quad T_{Y'}({\bf2_s})={3-\sqrt5\over4}
  \,,\quad T_{Y'}({\bf2'_s})={3+\sqrt5\over4}\cr
  T_{Y'}({\bf6_s})&=2\,,\quad T_{Y'}({\bf3'})=1-{1\over\sqrt5}
  \,,\quad T_{Y'}({\bf3})=1+{1\over\sqrt5}\cr
  T_{Y'}({\bf5})&=3/2,\quad T_{Y'}({\bf 4_s})=1\,.
  }
  \eql{twisted}
  $$
The new notation is that the rep is labelled by its dimension and, if a
spinor rep, has a suffix ${\bf s}$. Distinct reps with the same dimension
are distinguished by dashes and a dot means that the rep is one of a
conjugate pair, with the same torsion. \footnote{ We have adopted a
convention that the bundle is real if the twisting character is real, and
that a complex character is indicated by a power of two. An alternative
is to regard real twistings as particular complex ones, and then to
square every torsion. One could then elect {\it not} to square everything
and to work with an implied complex dimension. This we do later.}

It is interesting to note that the analytic torsions of Poincar\'e space,
for the spinor (quaternion), two--dimensional, double valued
representations, ${\bf2_s}$ and ${\bf2'_s}$, are half the conjugate
fundamental units of $\oQ(\sqrt5)$.\footnote{ The machine algebra we used
was not able to reduce to this form. We had to use Pell's equation.}
These values agree with the relations obtained by Tsuchiya,
[\pref{Tsuchiya}].

In order to amplify this, some explanatory remarks on signs and factors
of two are needed and the actual definition of torsion is relevant. We
have chosen the usual, but by no means universal, definition,
(\peq{at1}), employed by Ray and Singer, [\pref{RandS2}]. In terms of
this $T$, Ray's, [\pref{Ray}], definition is,
  $$
  T\big|_{Ray}=-2 \ln T=\ln(1/T^2)\,,
  \eql{tordefs}
  $$
assuming that the \zfs\ are the same. Since both references [\pref{Ray}]
and [\pref{RandS2}] deal with complex line bundles, this should be the
case. (See also [\pref{Rosenberg}].) However, on lens spaces, Ray
calculates, [\pref{Ray}], (1), (11),
  $$
  T_{Ray}=-2\sum_j \ln 2\sin( \pi r l_j /p)
  \eql{tray}
  $$
constructing the \zfs\ from the degeneracies of a {\it real} line bundle,
so far as we can see.

When Ray's result is referred to, it is usually the $T$ derived from
(\peq{tray}) using (\peq{tordefs}) that is  quoted. We note that this is
one half of our result, (\peq{tlens}), which incorporates a complex
dimension of two in the degeneracies; which is consistent.

In particular, Tsuchiya, [\pref{Tsuchiya}], states that the torsion on a
(one--sided) lens space, S$^3/\oZ_{10}$, is, from [\pref{Ray}],
  $$
  T_{Z_{10}}(\om^r)=|\om^r-1|^2\,,\quad \om^{10}=1\,,
  \eql{Ray3}
  $$
in conformity with (\peq{tordefs}) and (\peq{tray}), rather than with
(\peq{tlens})
\section{\bf 8. Clifford--Klein from lens by induction.}
The values of the torsion on Clifford--Klein spaces can be obtained from
those on lens spaces, although not {\it quite} for free, as claimed by
Ray, [\pref{Ray}]. The general result is guaranteed by Artin's theorem on
the rational sufficiency of the representations induced from all the
cyclic subgroups of the deck group and from the covering theorem of Ray
and Singer, [\pref{RandS2}], theorem 2.6 (see below).

Tsuchiya, [\pref{Tsuchiya}], on this basis, derives some relations
between torsions for the icosahedral case, which we summarise, and extend
here.

The basic connecting relation is the Ray--Singer covering theorem,
  $$
  T(\wt\man/\Ga_1;{\rm Ind}\, \rho)=T(\wt\man/\Ga_2;\rho)\,,\quad
  \Ga_2\subset\Ga_1\,,
  $$
applied, in particular to the cyclic subgroups of $\Ga$. For example, in
the present situation,
  $$
  \eqalign{
  T_{Y'}({\rm Ind}\, \rho)&=T_{Z_{10}}(\rho)\cr
  }
  \eql{thm26}
  $$
where $\oZ_{10},\,=\oZ_5\times\oZ_2,$ is generated by $A$ in the
Hamilton--Coxeter presentation of $Y'$, $(C^2=B^3= A^5, ABC)$.

Choosing an irrep for $\rho$, the dimension of the induced representation
is $120/10=12$, which is the size of the relevant classes in $Y'$. An
irrep, $\rho$, is generated by a power, $\om^r$, of a primitive tenth
root of unity, say $\om=\om_{10}=e^{\pi i/5}$, and a simple induced
character calculation yields,
  $$\eqalign{
  {\rm Ind}\, \om_{10}\,\,&=\,{\bf 2_s}\oplus{\bf 4_s}\oplus{\bf 6_s}\cr
  {\rm Ind}\, \om_{10}^3{\hspace{.}}&=\,{\bf 2'_s}\oplus{\bf 4_s}
  \oplus{\bf 6_s}\cr
  {\rm Ind}\, \om_{10}^5{\hspace{.}}&=\,{\bf 6_s}\oplus{\bf 6_s}\,.
  }
  $$

These decompositions translate into the torsion relations,
$$\eqalign{
  T_{Y'}({\bf 2_s})\,T_{Y'}({\bf 4_s})\,T_{Y'}({\bf 6_s})&
  =T_{Z_{10}}(\om_{10})=(3-\sqrt5)/2\cr
  T_{Y'}({\bf 2'_s})\,T_{Y'}({\bf 4_s})\,T_{Y'}({\bf 6_s})&
  =T_{Z_{10}}(\om_{10}^3)=(3+\sqrt5)/2\cr
  T_{Y'}({\bf 6_s})^2&=T_{Z_{10}}(\om_{10}^5)=4\,,
  }
  \eql{relns1}
  $$
the first two of which were given by Tsuchiya, except he does not
identify the ten--dimensional rep with ${\bf 4_s}\oplus{\bf 6_s}$.

The same calculation for $\oZ_6,\,=\oZ_3\times\oZ_2,$ generated by $B$,
produces the induction ${\rm Ind}\, \om_6={\bf 2_s}\oplus{\bf
2'_s}\oplus{\bf 4_s}\oplus{\bf 6_s}\oplus{\bf 6_s}$ and thence the
relation,
  $$\eqalign{
  T_{Y'}({\bf 2_s})\,T_{Y'}({\bf 2'_s})\,T_{Y'}({\bf 4_s})\,T^2_{Y'}({\bf 6_s})&
  =T_{Z_6}(\om_6)=1\,,\quad \om_6=e^{i\pi/3}\,.\cr
  }
  \eql{relns2}
  $$

All the spinor torsions in (\peq{twisted}) follow immediately from
(\peq{relns1}) and (\peq{relns2}). We note that, for numerical
consistency, the implied complex dimension convention is used.

An even $\oZ_{10}$ induction involves the non--spinor representations,
  $$
  {\rm Ind}\,\om_{10}^6={\bf 3}\oplus{\bf 4}\oplus{\bf 5}\,,\quad
  \om_{10}^6=-\om_{10}\,,
  $$
and it is left as an exercise to evaluate the non--spinor torsion values.

Nothing is gained by looking at say, $\oZ_5$, as this is a subgroup of
$\oZ_{10}$ and induction is transitive. Thus, for $\oZ_5$, the induced
reps are sums of $\oZ_{10}$ induced ones. For example ($\om_5=e^{2\pi
i/5}$),
  $$\eqalign{
  {\rm Ind}\, \om_5&={\rm Ind}\, \om_{10}\oplus{\rm Ind}
  ( -\om_{10})=
  \big({\bf 2_s}'\oplus{\bf 4_s}\oplus{\bf 6_s}\big)
  \oplus\big({\bf 3}\oplus{\bf 4}\oplus{\bf 5}\big)\cr
  {\rm Ind}\, \om_5^3&={\rm Ind}\, \om^3_{10}\oplus
  {\rm Ind}( -\om^3_{10})=\big({\bf 2_s{\hspace{'}}}
  \oplus{\bf 4_s}\oplus{\bf 6_s}\big)
  \oplus\big({\bf 3'}{\hspace{'}}\oplus{\bf 4}\oplus{\bf
  5}\big)\,,\cr
  }
  $$
and the equality of the analytic torsions is equivalent to the trivial
relation, $|\om^2-1|=|\om-1||-\om-1|$.

This induction method of deriving the torsion values is another, perhaps
more elegant, way of organising the information used in the brute force
evaluation via (\peq{torsion3}). Either way, one needs the character
tables.

A further check on the numbers is afforded by inducing from the {\it
trivial} representation. Any cyclic group would do; we select $\oZ_{10}$.
Then characters produce the decomposition,
  $$
  {\rm Ind}({\bf 1})={\bf 1}\oplus{\bf 3}\oplus{\bf 3'}\oplus{\bf 5}\,.
  \eql{trivind}
  $$
The lens space trivial torsion is given by (\peq{trivlens}), (not
computed by Ray),
  $$
  T_{Z_{10}}({\bf 1})={1\over100}\,|S^3|
  $$
whence, using (\peq{thm26}) and the values (\peq{twisted}), one finds
  $$
  T_{Y'}({\bf 1})={1\over120}\,|S^3|
  $$
agreeing with our earlier direct evaluation, (\peq{CK2}).

In the light of (\peq{relns1}), the question now arises of showing that
$T_{Z_{10}}(\om^3_{10})$ is the fundamental unit of $\oQ(\sqrt5)$, {\it
without} numerical evaluation. We do not pursue this point except to say
that the appearance of $(3\pm\sqrt5)/2$ is perhaps not unexpected in the
light of the general similarity of the torsion formula (\peq{torsion3})
to that for the class number of quadratic forms with positive
discriminant (\eg\ Lerch, [\pref{Lerch2}], p.366, Zagier,
[\pref{Zagierzf}], p.81).

\section {\bf 9. Polytope group averages.}
In general, the group average cannot be given in abstract closed form and
one has to resort to numerical addition class by class, at least in this
approach. This will be so for the finite (rotational) polytope groups,
$\{3,3,3\}$, $\{3,3,4\}$, $\{3,4,3\}$ and $\{3,3,5\}$.

The class decompositions for these have been given in [\pref{Chang}]. The
character tables of the full, reflective polytope groups, $[p,q,r]$, were
computed using the CAYLEY computer algebra system (this was superseded by
MAGMA in 1993) from their Coxeter presentations and the required subgroup
of index two selected. Much of the information can be found in Hurley,
[\pref{Hurley}], except for $\{3,3,5\}$, which is not crystallographic.
These character tables can also be found in Warner, [\pref{Warner}] and
some in Littlewood, [\pref{Littlewood}]. The relevant data have been
published in [\pref{Dow40}] and, being reasonably extensive, will not all
be repeated here and we can be brief. We will, however, demonstrate how
to calculate the torsion for $\{3,3,3\}$, the class structure of which is
  $$
   \{3,3,3\}=I\oplus15E\oplus 20K\oplus (2\times12)L'\,.
  $$

The $b$ coefficients\footnote{ The $b$s are not degeneracies.} in
(\peq{torsion3}), calculated using the angles $\al$ and $\be$ for the
non--trivial classes, $E,\,K$ and $L'$ are $(0,4)$, $(1,1,4)$ and
$(-1,-1,-1,-1,4)$ respectively, where we have used the order, $q$, for
each class, of $2$, $3$ and $5$. Putting the ingredients together we
find,
  $$\eqalign{
  \ln T_{\{3,3,3\}}(\rho)&={1\over60}\bigg[\ln2\pi^2\sum_\ga\chi^*(\ga)
  -30\chi_E^*\ln2-20\chi_K^*
  \big(T_3(1)+T_3(2)+2\ln3\big)\cr
  &12(\chi^*_{L_1'}+\chi^*_{L_2'})\big(T_5(1)+T_5(2)+T_5(3)+T_5(4)-2\ln5\big)
  \bigg]\cr
  &={1\over60}\ln2\pi^2\sum_\ga\chi^*(\ga)-{1\over2}\big(\chi_E^*\ln2+\chi_K^*
  \ln3  +(\chi^*_{L_1'}+\chi^*_{L_2'}\big)\ln5\big)\,,
  }
  \eql{333}
  $$
where we have used the product, (\peq{oldprod}), to simplify things. The
twisting, $\chi$, is chosen from the character tables.

The first term on the right of (\peq{333}) vanishes for a non--trivial
twisting, of which there are four irreducible ones. The corresponding
values of the torsion $T$ are $\sqrt{2/5}$ for both three dimensional
irreps, $5/\sqrt3$ and $\sqrt{3/2}$ for the four--dimensional and
five--dimensional ones, respectively.

We list the analytic torsions for all the twelve irreps of the
$\{3,3,4\}$ group. $\sqrt{8/3}$, $\sqrt{3/2}$, $1$, $1/\sqrt2$, $1$,
$1/\sqrt2$, $\sqrt8$, $1/\sqrt4$, $\sqrt4$, $\sqrt2$, $\sqrt2$,
$1/\sqrt2$.

For the {\it trivial} representation, one finds that, for the doubled
orbifold fundamental domains,
  $$\eqalign{
   \ln T_{\{3,3,3\}}({\bf 1})&=\ln
   (2\pi^2/60)+{3\over2}\ln2+{1\over2}\ln3\cr
   \ln T_{\{3,3,4\}}({\bf 1})&=\ln
   (2\pi^2/192)+2\ln2+{1\over2}\ln3\cr
   }
   \eql{tor333}
  $$
the torsion parts of which are not of Cheeger's combinatorial form,
(\peq{flat2}). That there should be some difference is not unreasonable
in view of the calculations of L\"uck, [\pref{Luck}], and Lott and
Rothenberg, [\pref{LandR}] on the effect of boundaries and non--free
actions. We conjecture that any extra terms are consequences of the
codimension--two conical singularities. The absence of a $\ln5$ term is
noted for $\{3,3,3\}$.

\section{\bf 10. Conclusion.}

The explicit computations performed here show that the twisted torsion is
an algebraic number, as expected on general grounds.

On the doubled tessellations terms occur which we attribute to
codimension--two singularities.

An alternative technique, which avoids the need to sum over classes and
angles, is to describe the polytope groups in terms of their integer {\it
degrees}, as in [\pref{Dow40,dowtess2}] and use a twisted generating
function. This is left for another time.

It has to be said that our discussion is geared to the three--sphere
while Ray's applies to any dimension. The extension of our cancellation
to higher dimensions via a generating function is also left for another
time.

\section{\bf Appendix 1. The cancellation.}

We wish to make some technical remarks concerning the important
cancellation (\peq{canc}). The derivation given earlier is somewhat
cavalier and needs tightening. For example, the trigonometric factor in
the definition (\peq{eff}) which, for $\th_L\ne0$ or $\th_R\ne0$,
provides a regulation necessary for the absence of a determinant
multiplicative anomaly associated with the product of `eigenvalues',
$n^2-1=(n+1)(n-1)$. We flesh out this comment with some algebra.

There are several cases that should be considered separately. First,
assume that $\th_L\ne0$ {\it and} $\th_R\ne0$ \ie $\al\ne\pm\be$. Then
write,
  $$
  \chi_n(\th_L)\,\chi_n(\th_R)={\cos n\al-\cos
  n\be\over\cos\al-\cos\be}\,,
  $$
so that we can concentrate on just the factor $\cos n\al$.

To facilitate the analysis, we make use of the Lerch--Lipshitz \zf,
$\Phi$, or rather a function simply related to it. Thus we define,
  $$\eqalign{
  \Xi(\al,s,w)&=\sum_{n=2}^\infty {e^{in\al}\over(n+w)^s}=e^{i\al}
  \sum_{n=1}^\infty
  {e^{in\al}\over(n+w+1)^s}\cr
  &=e^{i\al}\bigg[\Phi(e^{i\al},s,w+1)-{1\over(w+1)^s}\bigg]\,.
  }
  \eql{xi}
  $$

The sums are defined only for $\Re s>0$ but $\Phi$ can be continued in a
known fashion (\eg\ [\pref{Lerch}], [\pref{EMOT}]).

For the third sum in $F$, we employ a very old technique, explained in
[{\pref{Dow3}], and more generally in [\pref{DoandKi}], which is
sufficient to continue to $s=0$, at least.

To this end, we write it as,
  $$
  -2\lim_{u\to1}\Re\sum_2^\infty{e^{in\al}\over (n^2-u^2)^s}\equiv -2\lim_{u\to1}
  \Re\,\Upsilon(\al,2s,u)\,,\quad 0\le |u|<2\,.
  $$

In terms of $\Xi$ and $\Upsilon$, the essential part of $F$, (\peq{eff}),
is
  $$
  F(\al,s)=\Re\lim_{u\to1}\big(\Xi(\al,2s,u)
  +\Xi(\al,2s,-u)-2\Upsilon(\al,2s,u)\big)\,.
  \eql{ftwo2}
  $$

The $\Re$ is only a convenient bookkeeping symbol, not applying to the
complex nature of $s$. It will sometimes be removed and replaced at the
end of manipulations.

Expanding $\Upsilon$ in powers of $u^2$ one gets,
  $$
  \Upsilon(\al,2s,u)=\sum_{r=0}^\infty
  u^{2r}{s(s+1)\ldots(s+r-1)\over
  r!}\,\Upsilon(\al,2s+2r,0)
  \eql{th}
  $$
where $\Upsilon(*,*,0)$ is just $\Xi(*,*,0)$,
  $$
  \Upsilon(\al,2s,0)=\Xi(\al,2s,0)\,.
  \eql{up0}
  $$
This yields values at $s=0$ because the continuation of $\Xi$ is known.
For $\al\ne0$ it is an entire function of $s$ and can be given as a
contour integral, [\pref{Lerch}].

For example, at $s=0$, the absence of poles in $\Xi$ means that,
  $$
  -2\Upsilon(\al,0,u)=-2\Xi(\al,0,0)\,,\quad 0\le|u|<2\,,
  $$
and actual calculation shows, correctly, that this cancels the first two
sums in (\peq{eff}), \ie (\peq{ftwo2}). We soon show that such an {\it
explicit} demonstration is unnecessary.

Next the derivative at $s=0$ is required, and, again, the absence of
poles in $\Xi$, implies,
  $$
  2\Upsilon'(\al,0,u)=2\Xi'(\al,0,0)+\sum_{r=1}^\infty
  {u^{2r}\over r}\,\Xi(\al,2r,0)\,.
  \eql{th2}
  $$
using (\peq{up0}).

To evaluate the final summation in (\peq{th2}), we proceed as in
[\pref{Dow3}] and write $\Xi$ as an integral of the corresponding
cylinder kernel, $T$,
  $$\eqalign{
  \sum_{r=1}^\infty {u^{2r}\over r}\,\Xi(\al,2r,0)&=
  \sum_{r=1}^\infty{u^{2r}\over r\Ga(2r)}\int_0^\infty dt
  t^{2r-1}e^{-t}T(t,\al)\cr
  &=\int_0^\infty(\cosh ut-1)e^{-t}T(t,\al)\,{dt\over t}\cr
  &={1\over2}\lim_{s\to0}\int_0^\infty(e^{ut}
  +e^{-ut}-2)t^{s-1}\,e^{-t}T(t,\al)\,dt\,,
  }
  \eql{last}
  $$
where $T$ is contained in Lipshitz' formula, [\pref{EMOT}] 1.11 (4)
  $$
  \Xi(\al,s,w)={1\over2\Ga(s)}\int_0^\infty t^{s-1}e^{-wt}e^{-t}
  {e^{i\al}(e^{i\al}-e^{-t})\over\cosh t-\cos\al}\,,
  $$
and an auxiliary regularisation has been put in to assist the evaluation
of the (finite) integral.

A further, reverse, application  of Lipshitz' formula yields, simply
  $$
  \sum_{r=1}^\infty {u^{2r}\over
  r}\,\Xi(\al,2r,0)=\lim_{s\to0}
  \Ga(s)\big(\Xi(\al,s,u)+\Xi(\al,s,-u)-2\Xi(\al,s,0)\big)\,.
  \eql{last2}
  $$

This quantity must be finite, and so, in particular,
  $$
  \Xi(\al,0,1)+\Xi(\al,0,-1)-2\,\Xi(\al,0,0)=0
  $$
which expresses the vanishing of $F(s)$ at $s=0$ without explicit
evaluation, as advertised above.

The value of (\peq{last2}) is, therefore,
  $$
  \sum_{r=1}^\infty {u^{2r}\over
  r}\,\Xi(\al,2r,0)=\Xi'(\al,0,u)+\Xi'(\al,0,-u)-2\Xi'(\al,0,0)
  $$
and simple algebra shows that this is equivalent to
  $$
  F'(\al,0)=0
  $$
from (\peq{ftwo2}), as required.

If $\al=0$, only a slight modification is needed. Then $\Xi(0,s,w)$ is
essentially the Hurwitz \zf\ and has a single pole, at $s=1$. However
this does not change anything since only poles at positive even integers
can contribute to the evaluation of (\peq{th}) at $s=0$. The
multiplicative anomaly still vanishes.

This derivation of the cancellation is not much different from the rather
crude one in the main body of this paper.

The other cases should be examined separately. If $\al=\pm\be\ne0$, \ie
the group action is either all left or all right, then
  $$
  \chi_n(\th_L)\,\chi_n(\th_R)=n{\sin
  n\al\over\sin\al}=-{1\over\sin\al}{d\cos n\al\over d\al}
  $$
and the previous analysis can be applied since it is valid for all $\al$.

The same argument holds when $\th_L=\th_R=0$, which gives just the unit
element contribution. The degeneracies are the usual polynomials but all
one needs is the statement,
  $$
  \chi_n(0)\chi_n(0)=n^2=-{d^2\over d\al^2}\cos n\al\bigg|_{\al=0}\,.
  $$

To convince ourselves that this is true, it is easy to perform an
explicit calculation. The quantity involved is,
  $$\eqalign{
  \sum_{n=2}^\infty{n^2\over(n-1)^{2s}}&+{n^2\over(n+1)^{2s}}-2{n^2\over
  (n^2-1)^s}\cr
  &=2\sum_{n=1}^\infty\bigg({n^2+1\over n^{2s}}
  -{n^2\over(n^2-1)^{s}}\bigg)-{1\over 2^{2s}}\cr
  &=2\,\ze_R(2s-2)+2\ze_R(2s)-{1\over 2^{2s}}-2\ze_{3}(s)\,,
  }
  \eql{full}
  $$
where $\ze_{3}(s)$ is the minimal scalar \zf\ on the full three sphere.

The derivative of this at $0$ is
  $$
  4\ze_R'(-2)+4\ze_R'(0)+2\ln2-2\ze_3'(0)\,,
  $$
which vanishes using the known value of $\ze_3'(0)$ given in
[\pref{Dowjmp}], for example.\footnote{ This is really no  more than an
algebraic check as the method in [\pref{Dowjmp}] is just that used here.
However there are other techniques, such as Plana summation, which yield
the same value.}

\section{\bf Appendix 2. Cyclic decompositions.}

Due to the cancellation, the effective eigenvalues are squares of
integers, as in (\peq{aux}). This means that the calculations in
[\pref{ChandD,Dow11}] on spectral problems with conformal operators in
spherical factors can be applied directly to the present situation.

It was shown in [\pref{Dow11}] that the \zfs\ on fixed point free, one
sided factors, S$^3/\Ga'$, \ie homogeneous Clifford--Klein spaces, are
related to the zeta functions on cyclic factors (one sided lens spaces)
by going through an S$^2$ orbifold intermediary. Applying this result to
the (total) effective tau function gives,
  $$
  \wt\tau_{\Ga'}(s)={1\over2}
  \bigg(\sum_q\wt\tau_{Z_{2q}}(s)-\wt\tau_{Z_2}(s)\bigg)
  \eql{clkl}
  $$
where $\Ga'$ is a binary polyhedral group, one of $T'$, $O'$ or $ Y'$.
The sum is over the set of conjugate $q$--fold rotation axes appropriate
to the SO(3) {\it rotation} groups, ${\bf T}$, ${\bf O}$ or ${\bf Y}$.
(See [\pref{ChandD}].) The half and the fact that only even dimensional
lens spaces occur are consequences of the binary doubling. The
orbit--stabiliser relation has been used in reaching (\peq{clkl}).

For technical reasons, only untwisted fields were considered in
[\pref{Dow11}] and so the result (\peq{clkl}) applies immediately only to
trivial bundles, for which we have the lens space values
(\peq{trivlens}). The values of $q$ are, in each case, ${\bf T}:2,3,3$,
${\bf O}:2,3,4$ and ${\bf Y}: 2,3,5$. The analytic torsion is combined in
the same way as (\peq{clkl}). For example,
  $$\eqalign{
  \ln T_{T'}&={1\over2}\bigg(2\ln|S^3|-2\ln (4.6.6)+2\ln2\bigg)\cr
  &=\ln|S^3|-\ln (4.4.6)+\ln2\cr
  &=\ln\big(|S^3|/24\big)+\ln24-\ln (4.6.6)+\ln2\cr
  &=\ln\big(|S^3|/24\big)-\ln3\,,
  }
  $$
and similarly,
  $$\eqalign{
  \ln T_{O'}&=\ln\big(|S^3|/48\big)-\ln2\,,\quad
  \ln T_{Y'}=\ln\big(|S^3|/120\big)\,,\cr
  }
  $$
in agreement with our previous results, (\peq{CK2}).

The orbit--stabiliser relation does not hold for the cyclic case itself,
nor for the dihedral group. For the latter, one has the relation for the
conformal \zfs,\footnote{ One could just set $q=2,2,q$ in (\peq{clkl}).}
  $$
 \ze_{D'_{2q}}={1\over2}\big(\ze_{Z_{2q}}+2\ze_{Z_4}-\ze_{Z_2}\big)\,,
  $$
 and likewise for the effective $\tau$ function. Calculation rapidly
 reproduces (\peq{dihedral}).
 \section{\bf Appendix 3. Some homology.}
The relation obtained in [\pref{Dow11}] between the conformal \zfs\ on
one--sided factors S$^3/\Ga'$ and on the orbifolds, S$^2/\Ga$, is
parallelled by a relation between the homologies of fundamental domains.
Here we just discuss the first homology group, which is sufficient for
three dimensions.

The natural action of the rotation group, $\Ga$, on $\oR^3$ has a
(double) M\"obius corner as its (infinite) fundamental domain,
$\man=\oR^3/\Ga$, (\cf [\pref{DandCh}]). Denoting the generators of $\Ga$
by $A_i$ for rotations through the angles $2\pi/\nu_i$, about the
vertices of a given spherical triangle on S$^2$, we consider the
presentation $\Ga=(A_i:A_i^{\nu_i},A_1A_2A_3)$. Thus
$H_1(\man.\oZ)\cong\Ga/F$ where $F$ is the commutator subgroup of $\Ga$,
that is
  $$
  H_1(\man,\oZ)=\big(A_i:A_i^{\nu_i},A_1A_2A_3,
  A_1A_2A_1^{-1}A_2^{-1},(1\!\to\!2,2\!\to\!3),(1\!\to3,2\!\to1)\big)
  $$
from which, one finds
  $$
  H_1(\oR^3/\Ga,\oZ)=\oZ_3,\oZ_2, \{{\rm id}\}\quad{\rm for}\,\,{\bf
  T,O,Y}\quad {\rm respectively}\,.
  \eql{hommob}
  $$

The cyclic group, $\oZ_q$ is already abelian and so, elementarily,
  $$
  H_1(\oR^3/\oZ_q,\oZ)=\oZ_q\,.
  $$

These results coincide with those for the corresponding Clifford--Klein
S$^3$ factors, that is, generalising slightly,
  $$
  H_*(\oR^3/\Ga,\oZ)=H_*(S^3/\Ga',\oZ)\,.
  $$

A U(1) bundle (a complex field) can be twisted by a representation of
$H_1$ through phase factors around the triangle vertices,
 $$
  a(\ga)=e^{2\pi i{\bf n.\bmu}}
 $$
where ${\bf n}=(n_1,n_2,n_3)$ is the set of occurrences of the
generators, $A_1,A_2,A_3$, in the word presentation of $\ga$ and
$\mu_i=s_i/\nu_i$.

The relation condition $a(A_1A_2A_3)=a(E)=1$ translates into $\sum_i
\mu_i\in\oZ$ with $0\le\mu_i<1$ which determines the possible sets ${\bf
s}=(s_1,s_2,s_3)$ for given $\nu_i$. One straightforwardly
finds,\footnote{ This just confirms the known character tables.} in
addition to the trivial solution $(0,0,0)$, the two (complex)
representations ${\bf s}=(0,2,1)$ and $(0,1,2)$ for ${\bf T}$ and the
real one, $(1,0,2)$, for ${\bf O}$. There are no such reps for ${\bf Y}$.
This result is equivalent to the homology (\peq{hommob}).

\newpage
\noin{\bf References.} \vskip5truept
\begin{putreferences}
   \ref{Zagierzf}{Zagier,D. {\it Zetafunktionen und Quadratische
   K\"orper}, (Springer--Verlag, Berlin, 1981).}
   \ref{BandF}{Bauer,W. and Furutani, K. {\it J.Geom. and Phys.} {\bf
   58} (2008) 64.}
   \ref{Luck}{L\"uck,W. \jdg{37}{1993}{263}.}
   \ref{LandR}{Lott,J. and Rothenberg,M. \jdg{34}{1991}{431}.}
   \ref{DoandKi} {Dowker.J.S. and Kirsten, K. {\it Analysis and Appl.}
   {\bf 3} (2005) 45.}
   \ref{dowtess1}{Dowker,J.S. \cqg{23}{2006}{1}.}
   \ref{dowtess2}{Dowker,J.S. {\it J.Geom. and Phys.} {\bf 57} (2007) 1505.}
   \ref{MHS}{De Melo,T., Hartmann,L. and Spreafico,M. {\it Reidemeister
   Torsion and analytic torsion of discs}, ArXiv:0811.3196.}
   \ref{Vertman}{Vertman, B. {\it Analytic Torsion of a  bounded
   generalized cone}, ArXiv:0808.0449.}
   \ref{WandY} {Weng,L. and You,Y., {\it Int.J. of Math.}{\bf 7} (1996)
   109.}
   \ref{ScandT}{Schwartz, A.S. and Tyupkin,Yu.S. \np{242}{1984}{436}.}
   \ref{AAR}{Andrews, G.E., Askey,R. and Roy,R. {\it Special functions}
  (CUP, Cambridge, 1999).}
   \ref{Tsuchiya}{Tsuchiya, N. {\it Nagoya Math.J.} (1975) 289. }
  \ref{Lerch}{Lerch,M. \am{11}{1887}{19}.}
  \ref{Lerch2}{Lerch,M. \am{29}{1905}{333}.}
  \ref{TandS}{Threlfall, W. and Seifert, H. \ma{104}{1930}{1}.}
  \ref{RandS}{Ray, D.B., and Singer, I. \aim{7}{1971}{145}.}
  \ref{RandS2}{Ray, D.B., and Singer, I. {\it Proc.Symp.Pure Math.}
  {\bf 23} (1973) 167.}
  \ref{Jensen}{Jensen,J.L.W.V. \aom{17}{1915-1916}{124}.}
  \ref{Rosenberg}{Rosenberg, S. {\it The Laplacian on a Riemannian Manifold}
  (CUP, Cambridge, 1997).}
  \ref{Nando2}{Nash, C. and O'Connor, D-J. {\it Int.J.Mod.Phys.}
  {\bf A10} (1995) 1779.}
  \ref{Fock}{Fock,V. \zfp{98}{1935}{145}.}
  \ref{Levy}{Levy,M. \prs {204}{1950}{145}.}
  \ref{Schwinger2}{Schwinger,J. \jmp{5}{1964}{1606}.}
  \ref{Muller}{M\"uller, \lnm{}{}{}.}
  \ref{VMK}{Varshalovich.}
  \ref{DandWo}{Dowker,J.S. and Wolski, A. \prA{46}{1992}{6417}.}
  \ref{Zeitlin1}{Zeitlin,V. {\it Physica D} {\bf 49} (1991).  }
  \ref{Zeitlin0}{Zeitlin,V. {\it Nonlinear World} Ed by
   V.Baryakhtar {\it et al},  Vol.I p.717,  (World Scientific, Singapore, 1989).}
  \ref{Zeitlin2}{Zeitlin,V. \prl{93}{2004}{264501}. }
  \ref{Zeitlin3}{Zeitlin,V. \pla{339}{2005}{316}. }
  \ref{Groenewold}{Groenewold, H.J. {\it Physica} {\bf 12} (1946) 405.}
  \ref{Cohen}{Cohen, L. \jmp{7}{1966}{781}.}
  \ref{AandW}{Argawal G.S. and Wolf, E. \prD{2}{1970}{2161,2187,2206}.}
  \ref{Jantzen}{Jantzen,R.T. \jmp{19}{1978}{1163}.}
  \ref{Moses2}{Moses,H.E. \aop{42}{1967}{343}.}
  \ref{Carmeli}{Carmeli,M. \jmp{9}{1968}{1987}.}
  \ref{SHS}{Siemans,M., Hancock,J. and Siminovitch,D. {\it Solid State
  Nuclear Magnetic Resonance} {\bf 31}(2007)35.}
 \ref{Dowk}{Dowker,J.S. \prD{28}{1983}{3013}.}
 \ref{Heine}{Heine, E. {\it Handbuch der Kugelfunctionen}
  (G.Reimer, Berlin. 1878, 1881).}
  \ref{Pockels}{Pockels, F. {\it \"Uber die Differentialgleichung $\De
  u+k^2u=0$} (Teubner, Leipzig. 1891).}
  \ref{Hamermesh}{Hamermesh, M., {\it Group Theory} (Addison--Wesley,
  Reading. 1962).}
  \ref{Racah}{Racah, G. {\it Group Theory and Spectroscopy}
  (Princeton Lecture Notes, 1951). }
  \ref{Gourdin}{Gourdin, M. {\it Basics of Lie Groups} (Editions
  Fronti\'eres, Gif sur Yvette. 1982.)}
  \ref{Clifford}{Clifford, W.K. \plms{2}{1866}{116}.}
  \ref{Story2}{Story, W.E. \plms{23}{1892}{265}.}
  \ref{Story}{Story, W.E. \ma{41}{1893}{469}.}
  \ref{Poole}{Poole, E.G.C. \plms{33}{1932}{435}.}
  \ref{Dickson}{Dickson, L.E. {\it Algebraic Invariants} (Wiley, N.Y.
  1915).}
  \ref{Dickson2}{Dickson, L.E. {\it Modern Algebraic Theories}
  (Sanborn and Co., Boston. 1926).}
  \ref{Hilbert2}{Hilbert, D. {\it Theory of algebraic invariants} (C.U.P.,
  Cambridge. 1993).}
  \ref{Olver}{Olver, P.J. {\it Classical Invariant Theory} (C.U.P., Cambridge.
  1999.)}
  \ref{AST}{A\v{s}erova, R.M., Smirnov, J.F. and Tolsto\v{i}, V.N. {\it
  Teoret. Mat. Fyz.} {\bf 8} (1971) 255.}
  \ref{AandS}{A\v{s}erova, R.M., Smirnov, J.F. \np{4}{1968}{399}.}
  \ref{Shapiro}{Shapiro, J. \jmp{6}{1965}{1680}.}
  \ref{Shapiro2}{Shapiro, J.Y. \jmp{14}{1973}{1262}.}
  \ref{NandS}{Noz, M.E. and Shapiro, J.Y. \np{51}{1973}{309}.}
  \ref{Cayley2}{Cayley, A. {\it Phil. Trans. Roy. Soc. Lond.}
  {\bf 144} (1854) 244.}
  \ref{Cayley3}{Cayley, A. {\it Phil. Trans. Roy. Soc. Lond.}
  {\bf 146} (1856) 101.}
  \ref{Wigner}{Wigner, E.P. {\it Gruppentheorie} (Vieweg, Braunschweig. 1931).}
  \ref{Sharp}{Sharp, R.T. \ajop{28}{1960}{116}.}
  \ref{Laporte}{Laporte, O. {\it Z. f. Naturf.} {\bf 3a} (1948) 447.}
  \ref{Lowdin}{L\"owdin, P-O. \rmp{36}{1964}{966}.}
  \ref{Ansari}{Ansari, S.M.R. {\it Fort. d. Phys.} {\bf 15} (1967) 707.}
  \ref{SSJR}{Samal, P.K., Saha, R., Jain, P. and Ralston, J.P. {\it
  Testing Isotropy of Cosmic Microwave Background Radiation},
  astro-ph/0708.2816.}
  \ref{Lachieze}{Lachi\'eze-Rey, M. {\it Harmonic projection and
  multipole Vectors}. astro- \break ph/0409081.}
  \ref{CHS}{Copi, C.J., Huterer, D. and Starkman, G.D.
  \prD{70}{2003}{043515}.}
  \ref{Jaric}{Jari\'c, J.P. {\it Int. J. Eng. Sci.} {\bf 41} (2003) 2123.}
  \ref{RandD}{Roche, J.A. and Dowker, J.S. \jpa{1}{1968}{527}.}
  \ref{KandW}{Katz, G. and Weeks, J.R. \prD{70}{2004}{063527}.}
  \ref{Waerden}{van der Waerden, B.L. {\it Die Gruppen-theoretische
  Methode in der Quantenmechanik} (Springer, Berlin. 1932).}
  \ref{EMOT}{Erdelyi, A., Magnus, W., Oberhettinger, F. and Tricomi, F.G. {
  \it Higher Transcendental Functions} Vol.1 (McGraw-Hill, N.Y. 1953).}
  \ref{Dowzilch}{Dowker, J.S. {\it Proc. Phys. Soc.} {\bf 91} (1967) 28.}
  \ref{DandD}{Dowker, J.S. and Dowker, Y.P. {\it Proc. Phys. Soc.}
  {\bf 87} (1966) 65.}
  \ref{DandD2}{Dowker, J.S. and Dowker, Y.P. \prs{}{}{}.}
  \ref{Dowk3}{Dowker,J.S. \cqg{7}{1990}{1241}.}
  \ref{Dowk5}{Dowker,J.S. \cqg{7}{1990}{2353}.}
  \ref{CoandH}{Courant, R. and Hilbert, D. {\it Methoden der
  Mathematischen Physik} vol.1 \break (Springer, Berlin. 1931).}
  \ref{Applequist}{Applequist, J. \jpa{22}{1989}{4303}.}
  \ref{Torruella}{Torruella, \jmp{16}{1975}{1637}.}
  \ref{Weinberg}{Weinberg, S.W. \pr{133}{1964}{B1318}.}
  \ref{Meyerw}{Meyer, W.F. {\it Apolarit\"at und rationale Curven}
  (Fues, T\"ubingen. 1883.) }
  \ref{Ostrowski}{Ostrowski, A. {\it Jahrsb. Deutsch. Math. Verein.} {\bf
  33} (1923) 245.}
  \ref{Kramers}{Kramers, H.A. {\it Grundlagen der Quantenmechanik}, (Akad.
  Verlag., Leipzig, 1938).}
  \ref{ZandZ}{Zou, W.-N. and Zheng, Q.-S. \prs{459}{2003}{527}.}
  \ref{Weeks1}{Weeks, J.R. {\it Maxwell's multipole vectors
  and the CMB}.  astro-ph/0412231.}
  \ref{Corson}{Corson, E.M. {\it Tensors, Spinors and Relativistic Wave
  Equations} (Blackie, London. 1950).}
  \ref{Rosanes}{Rosanes, J. \jram{76}{1873}{312}.}
  \ref{Salmon}{Salmon, G. {\it Lessons Introductory to the Modern Higher
  Algebra} 3rd. edn. \break (Hodges,  Dublin. 1876.)}
  \ref{Milnew}{Milne, W.P. {\it Homogeneous Coordinates} (Arnold. London. 1910).}
  \ref{Niven}{Niven, W.D. {\it Phil. Trans. Roy. Soc.} {\bf 170} (1879) 393.}
  \ref{Scott}{Scott, C.A. {\it An Introductory Account of
  Certain Modern Ideas and Methods in Plane Analytical Geometry,}
  (MacMillan, N.Y. 1896).}
  \ref{Bargmann}{Bargmann, V. \rmp{34}{1962}{300}.}
  \ref{Maxwell}{Maxwell, J.C. {\it A Treatise on Electricity and
  Magnetism} 2nd. edn. (Clarendon Press, Oxford. 1882).}
  \ref{BandL}{Biedenharn, L.C. and Louck, J.D. {\it Angular Momentum in Quantum Physics}
  (Addison-Wesley, Reading. 1981).}
  \ref{Weylqm}{Weyl, H. {\it The Theory of Groups and Quantum Mechanics}
  (Methuen, London. 1931).}
  \ref{Robson}{Robson, A. {\it An Introduction to Analytical Geometry} Vol I
  (C.U.P., Cambridge. 1940.)}
  \ref{Sommerville}{Sommerville, D.M.Y. {\it Analytical Conics} 3rd. edn.
   (Bell. London. 1933).}
  \ref{Coolidge}{Coolidge, J.L. {\it A Treatise on Algebraic Plane Curves}
  (Clarendon Press, Oxford. 1931).}
  \ref{SandK}{Semple, G. and Kneebone. G.T. {\it Algebraic Projective
  Geometry} (Clarendon Press, Oxford. 1952).}
  \ref{AandC}{Abdesselam A., and Chipalkatti, J. {\it The Higher
  Transvectants are redundant}, arXiv:0801.1533 [math.AG] 2008.}
  \ref{Elliott}{Elliott, E.B. {\it The Algebra of Quantics} 2nd edn.
  (Clarendon Press, Oxford. 1913).}
  \ref{Elliott2}{Elliott, E.B. \qjpam{48}{1917}{372}.}
  \ref{Howe}{Howe, R. \tams{313}{1989}{539}.}
  \ref{Clebsch}{Clebsch, A. \jram{60}{1862}{343}.}
  \ref{Prasad}{Prasad, G. \ma{72}{1912}{136}.}
  \ref{Dougall}{Dougall, J. \pems{32}{1913}{30}.}
  \ref{Penrose}{Penrose, R. \aop{10}{1960}{171}.}
  \ref{Penrose2}{Penrose, R. \prs{273}{1965}{171}.}
  \ref{Burnside}{Burnside, W.S. \qjm{10}{1870}{211}. }
  \ref{Lindemann}{Lindemann, F. \ma{23} {1884}{111}.}
  \ref{Backus}{Backus, G. {\it Rev. Geophys. Space Phys.} {\bf 8} (1970) 633.}
  \ref{Baerheim}{Baerheim, R. {\it Q.J. Mech. appl. Math.} {\bf 51} (1998) 73.}
  \ref{Lense}{Lense, J. {\it Kugelfunktionen} (Akad.Verlag, Leipzig. 1950).}
  \ref{Littlewood}{Littlewood, D.E. \plms{50}{1948}{349}.}
  \ref{Fierz}{Fierz, M. {\it Helv. Phys. Acta} {\bf 12} (1938) 3.}
  \ref{Williams}{Williams, D.N. {\it Lectures in Theoretical Physics} Vol. VII,
  (Univ.Colorado Press, Boulder. 1965).}
  \ref{Dennis}{Dennis, M. \jpa{37}{2004}{9487}.}
  \ref{Pirani}{Pirani, F. {\it Brandeis Lecture Notes on
  General Relativity,} edited by S. Deser and K. Ford. (Brandeis, Mass. 1964).}
  \ref{Sturm}{Sturm, R. \jram{86}{1878}{116}.}
  \ref{Schlesinger}{Schlesinger, O. \ma{22}{1883}{521}.}
  \ref{Askwith}{Askwith, E.H. {\it Analytical Geometry of the Conic
  Sections} (A.\&C. Black, London. 1908).}
  \ref{Todd}{Todd, J.A. {\it Projective and Analytical Geometry}.
  (Pitman, London. 1946).}
  \ref{Glenn}{Glenn. O.E. {\it Theory of Invariants} (Ginn \& Co, N.Y. 1915).}
  \ref{DowkandG}{Dowker, J.S. and Goldstone, M. \prs{303}{1968}{381}.}
  \ref{Turnbull}{Turnbull, H.A. {\it The Theory of Determinants,
  Matrices and Invariants} 3rd. edn. (Dover, N.Y. 1960).}
  \ref{MacMillan}{MacMillan, W.D. {\it The Theory of the Potential}
  (McGraw-Hill, N.Y. 1930).}
   \ref{Hobson}{Hobson, E.W. {\it The Theory of Spherical and Ellipsoidal Harmonics}
   C.U.P., Cambridge. 1931).}
  \ref{Hobson1}{Hobson, E.W. \plms {24}{1892}{55}.}
  \ref{GandY}{Grace, J.H. and Young, A. {\it The Algebra of Invariants}
  (C.U.P., Cambridge, 1903).}
  \ref{FandR}{Fano, U. and Racah, G. {\it Irreducible Tensorial Sets}
  (Academic Press, N.Y. 1959).}
  \ref{TandT}{Thomson, W. and Tait, P.G. {\it Treatise on Natural Philosophy}
  (Clarendon Press, Oxford. 1867).}
  \ref{Brinkman}{Brinkman, H.C. {\it Applications of spinor invariants in
atomic physics}, North Holland, Amsterdam 1956.}
  \ref{Kramers1}{Kramers, H.A. {\it Proc. Roy. Soc. Amst.} {\bf 33} (1930) 953.}
  \ref{DandP2}{Dowker,J.S. and Pettengill,D.F. \jpa{7}{1974}{1527}}
  \ref{Dowk1}{Dowker,J.S. \jpa{}{}{45}.}
  \ref{Dowk2}{Dowker,J.S. \aop{71}{1972}{577}}
  \ref{DandA}{Dowker,J.S. and Apps, J.S. \cqg{15}{1998}{1121}.}
  \ref{Weil}{Weil,A., {\it Elliptic functions according to Eisenstein
  and Kronecker}, Springer, Berlin, 1976.}
  \ref{Ling}{Ling,C-H. {\it SIAM J.Math.Anal.} {\bf5} (1974) 551.}
  \ref{Ling2}{Ling,C-H. {\it J.Math.Anal.Appl.}(1988).}
 \ref{BMO}{Brevik,I., Milton,K.A. and Odintsov, S.D. \aop{302}{2002}{120}.}
 \ref{KandL}{Kutasov,D. and Larsen,F. {\it JHEP} 0101 (2001) 1.}
 \ref{KPS}{Klemm,D., Petkou,A.C. and Siopsis {\it Entropy
 bounds, monoticity properties and scaling in CFT's}. hep-th/0101076.}
 \ref{DandC}{Dowker,J.S. and Critchley,R. \prD{15}{1976}{1484}.}
 \ref{AandD}{Al'taie, M.B. and Dowker, J.S. \prD{18}{1978}{3557}.}
 \ref{Dow1}{Dowker,J.S. \prD{37}{1988}{558}.}
 \ref{Dow30}{Dowker,J.S. \prD{28}{1983}{3013}.}
 \ref{DandK}{Dowker,J.S. and Kennedy,G. \jpa{}{1978}{}.}
 \ref{Dow2}{Dowker,J.S. \cqg{1}{1984}{359}.}
 \ref{DandKi}{Dowker,J.S. and Kirsten, K. {\it Comm. in Anal. and Geom.
 }{\bf7} (1999) 641.}
 \ref{DandKe}{Dowker,J.S. and Kennedy,G.\jpa{11}{1978}{895}.}
 \ref{Gibbons}{Gibbons,G.W. \pl{60A}{1977}{385}.}
 \ref{Cardy}{Cardy,J.L. \np{366}{1991}{403}.}
 \ref{ChandD}{Chang,P. and Dowker,J.S. \np{395}{1993}{407}.}
 \ref{DandC2}{Dowker,J.S. and Critchley,R. \prD{13}{1976}{224}.}
 \ref{Camporesi}{Camporesi,R. \prp{196}{1990}{1}.}
 \ref{BandM}{Brown,L.S. and Maclay,G.J. \pr{184}{1969}{1272}.}
 \ref{CandD}{Candelas,P. and Dowker,J.S. \prD{19}{1979}{2902}.}
 \ref{Unwin1}{Unwin,S.D. Thesis. University of Manchester. 1979.}
 \ref{Unwin2}{Unwin,S.D. \jpa{13}{1980}{313}.}
 \ref{DandB}{Dowker,J.S.and Banach,R. \jpa{11}{1978}{2255}.}
 \ref{Obhukov}{Obhukov,Yu.N. \pl{109B}{1982}{195}.}
 \ref{Kennedy}{Kennedy,G. \prD{23}{1981}{2884}.}
 \ref{CandT}{Copeland,E. and Toms,D.J. \np {255}{1985}{201}.}
 \ref{ELV}{Elizalde,E., Lygren, M. and Vassilevich,
 D.V. \jmp {37}{1996}{3105}.}
 \ref{Malurkar}{Malurkar,S.L. {\it J.Ind.Math.Soc} {\bf16} (1925/26) 130.}
 \ref{Glaisher}{Glaisher,J.W.L. {\it Messenger of Math.} {\bf18}
(1889) 1.} \ref{Anderson}{Anderson,A. \prD{37}{1988}{536}.}
 \ref{CandA}{Cappelli,A. and D'Appollonio,\pl{487B}{2000}{87}.}
 \ref{Wot}{Wotzasek,C. \jpa{23}{1990}{1627}.}
 \ref{RandT}{Ravndal,F. and Tollesen,D. \prD{40}{1989}{4191}.}
 \ref{SandT}{Santos,F.C. and Tort,A.C. \pl{482B}{2000}{323}.}
 \ref{FandO}{Fukushima,K. and Ohta,K. {\it Physica} {\bf A299} (2001) 455.}
 \ref{GandP}{Gibbons,G.W. and Perry,M. \prs{358}{1978}{467}.}
 \ref{Dow4}{Dowker,J.S..}
  \ref{Rad}{Rademacher,H. {\it Topics in analytic number theory,}
Springer-Verlag,  Berlin,1973.}
  \ref{Halphen}{Halphen,G.-H. {\it Trait\'e des Fonctions Elliptiques},
  Vol 1, Gauthier-Villars, Paris, 1886.}
  \ref{CandW}{Cahn,R.S. and Wolf,J.A. {\it Comm.Mat.Helv.} {\bf 51}
  (1976) 1.}
  \ref{Berndt}{Berndt,B.C. \rmjm{7}{1977}{147}.}
  \ref{Hurwitz}{Hurwitz,A. \ma{18}{1881}{528}.}
  \ref{Hurwitz2}{Hurwitz,A. {\it Mathematische Werke} Vol.I. Basel,
  Birkhauser, 1932.}
  \ref{Berndt2}{Berndt,B.C. \jram{303/304}{1978}{332}.}
  \ref{RandA}{Rao,M.B. and Ayyar,M.V. \jims{15}{1923/24}{150}.}
  \ref{Hardy}{Hardy,G.H. \jlms{3}{1928}{238}.}
  \ref{TandM}{Tannery,J. and Molk,J. {\it Fonctions Elliptiques},
   Gauthier-Villars, Paris, 1893--1902.}
  \ref{schwarz}{Schwarz,H.-A. {\it Formeln und
  Lehrs\"atzen zum Gebrauche..},Springer 1893.(The first edition was 1885.)
  The French translation by Henri Pad\'e is {\it Formules et Propositions
  pour L'Emploi...},Gauthier-Villars, Paris, 1894}
  \ref{Hancock}{Hancock,H. {\it Theory of elliptic functions}, Vol I.
   Wiley, New York 1910.}
  \ref{watson}{Watson,G.N. \jlms{3}{1928}{216}.}
  \ref{MandO}{Magnus,W. and Oberhettinger,F. {\it Formeln und S\"atze},
  Springer-Verlag, Berlin 1948.}
  \ref{Klein}{Klein,F. {\it Lectures on the Icosohedron}
  (Methuen, London, 1913).}
  \ref{AandL}{Appell,P. and Lacour,E. {\it Fonctions Elliptiques},
  Gauthier-Villars,
  Paris, 1897.}
  \ref{HandC}{Hurwitz,A. and Courant,C. {\it Allgemeine Funktionentheorie},
  Springer,
  Berlin, 1922.}
  \ref{WandW}{Whittaker,E.T. and Watson,G.N. {\it Modern analysis},
  Cambridge 1927.}
  \ref{SandC}{Selberg,A. and Chowla,S. \jram{227}{1967}{86}. }
  \ref{zucker}{Zucker,I.J. {\it Math.Proc.Camb.Phil.Soc} {\bf 82 }(1977)
  111.}
  \ref{glasser}{Glasser,M.L. {\it Maths.of Comp.} {\bf 25} (1971) 533.}
  \ref{GandW}{Glasser, M.L. and Wood,V.E. {\it Maths of Comp.} {\bf 25}
  (1971)
  535.}
  \ref{greenhill}{Greenhill,A,G. {\it The Applications of Elliptic
  Functions}, MacMillan, London, 1892.}
  \ref{Weierstrass}{Weierstrass,K. {\it J.f.Mathematik (Crelle)}
{\bf 52} (1856) 346.}
  \ref{Weierstrass2}{Weierstrass,K. {\it Mathematische Werke} Vol.I,p.1,
  Mayer u. M\"uller, Berlin, 1894.}
  \ref{Fricke}{Fricke,R. {\it Die Elliptische Funktionen und Ihre Anwendungen},
    Teubner, Leipzig. 1915, 1922.}
  \ref{Konig}{K\"onigsberger,L. {\it Vorlesungen \"uber die Theorie der
 Elliptischen Funktionen},  \break Teubner, Leipzig, 1874.}
  \ref{Milne}{Milne,S.C. {\it The Ramanujan Journal} {\bf 6} (2002) 7-149.}
  \ref{Schlomilch}{Schl\"omilch,O. {\it Ber. Verh. K. Sachs. Gesell. Wiss.
  Leipzig}  {\bf 29} (1877) 101-105; {\it Compendium der h\"oheren
  Analysis}, Bd.II, 3rd Edn, Vieweg, Brunswick, 1878.}
  \ref{BandB}{Briot,C. and Bouquet,C. {\it Th\`eorie des Fonctions
  Elliptiques}, Gauthier-Villars, Paris, 1875.}
  \ref{Dumont}{Dumont,D. \aim {41}{1981}{1}.}
  \ref{Andre}{Andr\'e,D. {\it Ann.\'Ecole Normale Superior} {\bf 6} (1877)
  265;
  {\it J.Math.Pures et Appl.} {\bf 5} (1878) 31.}
  \ref{Raman}{Ramanujan,S. {\it Trans.Camb.Phil.Soc.} {\bf 22} (1916) 159;
 {\it Collected Papers}, Cambridge, 1927}
  \ref{Weber}{Weber,H.M. {\it Lehrbuch der Algebra} Bd.III, Vieweg,
  Brunswick 190  3.}
  \ref{Weber2}{Weber,H.M. {\it Elliptische Funktionen und algebraische
  Zahlen},
  Vieweg, Brunswick 1891.}
  \ref{ZandR}{Zucker,I.J. and Robertson,M.M.
  {\it Math.Proc.Camb.Phil.Soc} {\bf 95 }(1984) 5.}
  \ref{JandZ1}{Joyce,G.S. and Zucker,I.J.
  {\it Math.Proc.Camb.Phil.Soc} {\bf 109 }(1991) 257.}
  \ref{JandZ2}{Zucker,I.J. and Joyce.G.S.
  {\it Math.Proc.Camb.Phil.Soc} {\bf 131 }(2001) 309.}
  \ref{zucker2}{Zucker,I.J. {\it SIAM J.Math.Anal.} {\bf 10} (1979) 192,}
  \ref{BandZ}{Borwein,J.M. and Zucker,I.J. {\it IMA J.Math.Anal.} {\bf 12}
  (1992) 519.}
  \ref{Cox}{Cox,D.A. {\it Primes of the form $x^2+n\,y^2$}, Wiley,
  New York, 1989.}
  \ref{BandCh}{Berndt,B.C. and Chan,H.H. {\it Mathematika} {\bf42} (1995)
  278.}
  \ref{EandT}{Elizalde,R. and Tort.hep-th/}
  \ref{KandS}{Kiyek,K. and Schmidt,H. {\it Arch.Math.} {\bf 18} (1967) 438.}
  \ref{Oshima}{Oshima,K. \prD{46}{1992}{4765}.}
  \ref{greenhill2}{Greenhill,A.G. \plms{19} {1888} {301}.}
  \ref{Russell}{Russell,R. \plms{19} {1888} {91}.}
  \ref{BandB}{Borwein,J.M. and Borwein,P.B. {\it Pi and the AGM}, Wiley,
  New York, 1998.}
  \ref{Resnikoff}{Resnikoff,H.L. \tams{124}{1966}{334}.}
  \ref{vandp}{Van der Pol, B. {\it Indag.Math.} {\bf18} (1951) 261,272.}
  \ref{Rankin}{Rankin,R.A. {\it Modular forms} C.U.P. Cambridge}
  \ref{Rankin2}{Rankin,R.A. {\it Proc. Roy.Soc. Edin.} {\bf76 A} (1976) 107.}
  \ref{Skoruppa}{Skoruppa,N-P. {\it J.of Number Th.} {\bf43} (1993) 68 .}
  \ref{Down}{Dowker.J.S. \np {104}{2002}{153}.}
  \ref{Eichler}{Eichler,M. \mz {67}{1957}{267}.}
  \ref{Zagier}{Zagier,D. \invm{104}{1991}{449}.}
  \ref{Lang}{Lang,S. {\it Modular Forms}, Springer, Berlin, 1976.}
  \ref{Kosh}{Koshliakov,N.S. {\it Mess.of Math.} {\bf 58} (1928) 1.}
  \ref{BandH}{Bodendiek, R. and Halbritter,U. \amsh{38}{1972}{147}.}
  \ref{Smart}{Smart,L.R., \pgma{14}{1973}{1}.}
  \ref{Grosswald}{Grosswald,E. {\it Acta. Arith.} {\bf 21} (1972) 25.}
  \ref{Kata}{Katayama,K. {\it Acta Arith.} {\bf 22} (1973) 149.}
  \ref{Ogg}{Ogg,A. {\it Modular forms and Dirichlet series} (Benjamin,
  New York,
   1969).}
  \ref{Bol}{Bol,G. \amsh{16}{1949}{1}.}
  \ref{Epstein}{Epstein,P. \ma{56}{1903}{615}.}
  \ref{Petersson}{Petersson.}
  \ref{Serre}{Serre,J-P. {\it A Course in Arithmetic}, Springer,
  New York, 1973.}
  \ref{Schoenberg}{Schoenberg,B., {\it Elliptic Modular Functions},
  Springer, Berlin, 1974.}
  \ref{Apostol}{Apostol,T.M. \dmj {17}{1950}{147}.}
  \ref{Ogg2}{Ogg,A. {\it Lecture Notes in Math.} {\bf 320} (1973) 1.}
  \ref{Knopp}{Knopp,M.I. \dmj {45}{1978}{47}.}
  \ref{Knopp2}{Knopp,M.I. \invm {}{1994}{361}.}
  \ref{LandZ}{Lewis,J. and Zagier,D. \aom{153}{2001}{191}.}
  \ref{DandK1}{Dowker,J.S. and Kirsten,K. {\it Elliptic functions and
  temperature inversion symmetry on spheres} hep-th/.}
  \ref{HandK}{Husseini and Knopp.}
  \ref{Kober}{Kober,H. \mz{39}{1934-5}{609}.}
  \ref{HandL}{Hardy,G.H. and Littlewood, \am{41}{1917}{119}.}
  \ref{Watson}{Watson,G.N. \qjm{2}{1931}{300}.}
  \ref{SandC2}{Chowla,S. and Selberg,A. {\it Proc.Nat.Acad.} {\bf 35}
  (1949) 371.}
  \ref{Landau}{Landau, E. {\it Lehre von der Verteilung der Primzahlen},
  (Teubner, Leipzig, 1909).}
  \ref{Berndt4}{Berndt,B.C. \tams {146}{1969}{323}.}
  \ref{Berndt3}{Berndt,B.C. \tams {}{}{}.}
  \ref{Bochner}{Bochner,S. \aom{53}{1951}{332}.}
  \ref{Weil2}{Weil,A.\ma{168}{1967}{}.}
  \ref{CandN}{Chandrasekharan,K. and Narasimhan,R. \aom{74}{1961}{1}.}
  \ref{Rankin3}{Rankin,R.A. {} {} ().}
  \ref{Berndt6}{Berndt,B.C. {\it Trans.Edin.Math.Soc}.}
  \ref{Elizalde}{Elizalde,E. {\it Ten Physical Applications of Spectral
  Zeta Function Theory}, \break (Springer, Berlin, 1995).}
  \ref{Allen}{Allen,B., Folacci,A. and Gibbons,G.W. \pl{189}{1987}{304}.}
  \ref{Krazer}{Krazer}
  \ref{Elizalde3}{Elizalde,E. {\it J.Comp.and Appl. Math.} {\bf 118}
  (2000) 125.}
  \ref{Elizalde2}{Elizalde,E., Odintsov.S.D, Romeo, A. and Bytsenko,
  A.A and
  Zerbini,S.
  {\it Zeta function regularisation}, (World Scientific, Singapore,
  1994).}
  \ref{Eisenstein}{Eisenstein}
  \ref{Hecke}{Hecke,E. \ma{112}{1936}{664}.}
  \ref{Hecke2}{Hecke,E. \ma{112}{1918}{398}.}
  \ref{Terras}{Terras,A. {\it Harmonic analysis on Symmetric Spaces} (Springer,
  New York, 1985).}
  \ref{BandG}{Bateman,P.T. and Grosswald,E. {\it Acta Arith.} {\bf 9}
  (1964) 365.}
  \ref{Deuring}{Deuring,M. \aom{38}{1937}{585}.}
  \ref{Guinand}{Guinand.}
  \ref{Guinand2}{Guinand.}
  \ref{Minak}{Minakshisundaram.}
  \ref{Mordell}{Mordell,J. \prs{}{}{}.}
  \ref{GandZ}{Glasser,M.L. and Zucker, {}.}
  \ref{Landau2}{Landau,E. \jram{}{1903}{64}.}
  \ref{Kirsten1}{Kirsten,K. \jmp{35}{1994}{459}.}
  \ref{Sommer}{Sommer,J. {\it Vorlesungen \"uber Zahlentheorie}
  (1907,Teubner,Leipzig).
  French edition 1913 .}
  \ref{Reid}{Reid,L.W. {\it Theory of Algebraic Numbers},
  (1910,MacMillan,New York).}
  \ref{Milnor}{Milnor, J. {\it Is the Universe simply--connected?},
  IAS, Princeton, 1978.}
  \ref{Milnor2}{Milnor, J. \ajm{79}{1957}{623}.}
  \ref{Opechowski}{Opechowski,W. {\it Physica} {\bf 7} (1940) 552.}
  \ref{Bethe}{Bethe, H.A. \zfp{3}{1929}{133}.}
  \ref{LandL}{Landau, L.D. and Lishitz, E.M. {\it Quantum
  Mechanics} (Pergamon Press, London, 1958).}
  \ref{GPR}{Gibbons, G.W., Pope, C. and R\"omer, H., \np{157}{1979}{377}.}
  \ref{Jadhav}{Jadhav,S.P. PhD Thesis, University of Manchester 1990.}
  \ref{DandJ}{Dowker,J.S. and Jadhav, S. \prD{39}{1989}{1196}.}
  \ref{CandM}{Coxeter, H.S.M. and Moser, W.O.J. {\it Generators and
  relations of finite groups} Springer. Berlin. 1957.}
  \ref{Coxeter2}{Coxeter, H.S.M. {\it Regular Complex Polytopes},
   (Cambridge University Press,
  Cambridge, 1975).}
  \ref{Coxeter}{Coxeter, H.S.M. {\it Regular Polytopes}.}
  \ref{Stiefel}{Stiefel, E., J.Research NBS {\bf 48} (1952) 424.}
  \ref{BandS}{Brink, D.M. and Satchler, G.R. {\it Angular momentum theory}.
  (Clarendon Press, Oxford. 1962.).}
  \ref{Rose}{Rose}
  \ref{Schwinger}{Schwinger, J. {\it On Angular Momentum} in {\it Quantum Theory of
  Angular Momentum} edited by Biedenharn,L.C. and van Dam, H.
  (Academic Press, N.Y. 1965).}
  \ref{Bromwich}{Bromwich, T.J.I'A. {\it Infinite Series},
  (Macmillan, 1947).}
  \ref{Ray}{Ray,D.B. \aim{4}{1970}{109}.}
  \ref{Ikeda}{Ikeda,A. {\it Kodai Math.J.} {\bf 18} (1995) 57.}
  \ref{Kennedy}{Kennedy,G. \prD{23}{1981}{2884}.}
  \ref{Ellis}{Ellis,G.F.R. {\it General Relativity} {\bf2} (1971) 7.}
  \ref{Dow8}{Dowker,J.S. \cqg{20}{2003}{L105}.}
  \ref{IandY}{Ikeda, A and Yamamoto, Y. \ojm {16}{1979}{447}.}
  \ref{BandI}{Bander,M. and Itzykson,C. \rmp{18}{1966}{2}.}
  \ref{Schulman}{Schulman, L.S. \pr{176}{1968}{1558}.}
  \ref{Bar1}{B\"ar,C. {\it Arch.d.Math.}{\bf 59} (1992) 65.}
  \ref{Bar2}{B\"ar,C. {\it Geom. and Func. Anal.} {\bf 6} (1996) 899.}
  \ref{Vilenkin}{Vilenkin, N.J. {\it Special functions},
  (Am.Math.Soc., Providence, 1968).}
  \ref{Talman}{Talman, J.D. {\it Special functions} (Benjamin,N.Y.,1968).}
  \ref{Miller}{Miller, W. {\it Symmetry groups and their applications}
  (Wiley, N.Y., 1972).}
  \ref{Dow3}{Dowker,J.S. \cmp{162}{1994}{633}.}
  \ref{Cheeger}{Cheeger, J. \jdg {18}{1983}{575}.}
  \ref{Cheeger2}{Cheeger, J. \aom {109}{1979}{259}.}
  \ref{Dow6}{Dowker,J.S. \jmp{30}{1989}{770}.}
  \ref{Dow20}{Dowker,J.S. \jmp{35}{1994}{6076}.}
  \ref{Dowjmp}{Dowker,J.S. \jmp{35}{1994}{4989}.}
  \ref{Dow21}{Dowker,J.S. {\it Heat kernels and polytopes} in {\it
   Heat Kernel Techniques and Quantum Gravity}, ed. by S.A.Fulling,
   Discourses in Mathematics and its Applications, No.4, Dept.
   Maths., Texas A\&M University, College Station, Texas, 1995.}
  \ref{Dow9}{Dowker,J.S. \jmp{42}{2001}{1501}.}
  \ref{Dow7}{Dowker,J.S. \jpa{25}{1992}{2641}.}
  \ref{Warner}{Warner.N.P. \prs{383}{1982}{379}.}
  \ref{Wolf}{Wolf, J.A. {\it Spaces of constant curvature},
  (McGraw--Hill,N.Y., 1967).}
  \ref{Meyer}{Meyer,B. \cjm{6}{1954}{135}.}
  \ref{BandB}{B\'erard,P. and Besson,G. {\it Ann. Inst. Four.} {\bf 30}
  (1980) 237.}
  \ref{PandM}{Polya,G. and Meyer,B. \cras{228}{1948}{28}.}
  \ref{Springer}{Springer, T.A. Lecture Notes in Math. vol 585 (Springer,
  Berlin,1977).}
  \ref{SeandT}{Threlfall, H. and Seifert, W. \ma{104}{1930}{1}.}
  \ref{Hopf}{Hopf,H. \ma{95}{1925}{313}. }
  \ref{Dow}{Dowker,J.S. \jpa{5}{1972}{936}.}
  \ref{LLL}{Lehoucq,R., Lachi\'eze-Rey,M. and Luminet, J.--P. {\it
  Astron.Astrophys.} {\bf 313} (1996) 339.}
  \ref{LaandL}{Lachi\'eze-Rey,M. and Luminet, J.--P.
  \prp{254}{1995}{135}.}
  \ref{Schwarzschild}{Schwarzschild, K., {\it Vierteljahrschrift der
  Ast.Ges.} {\bf 35} (1900) 337.}
  \ref{Starkman}{Starkman,G.D. \cqg{15}{1998}{2529}.}
  \ref{LWUGL}{Lehoucq,R., Weeks,J.R., Uzan,J.P., Gausman, E. and
  Luminet, J.--P. \cqg{19}{2002}{4683}.}
  \ref{Dow10}{Dowker,J.S. \prD{28}{1983}{3013}.}
  \ref{BandD}{Banach, R. and Dowker, J.S. \jpa{12}{1979}{2527}.}
  \ref{Jadhav2}{Jadhav,S. \prD{43}{1991}{2656}.}
  \ref{Gilkey}{Gilkey,P.B. {\it Invariance theory,the heat equation and
  the Atiyah--Singer Index theorem} (CRC Press, Boca Raton, 1994).}
  \ref{BandY}{Berndt,B.C. and Yeap,B.P. {\it Adv. Appl. Math.}
  {\bf29} (2002) 358.}
  \ref{HandR}{Hanson,A.J. and R\"omer,H. \pl{80B}{1978}{58}.}
  \ref{Hill}{Hill,M.J.M. {\it Trans.Camb.Phil.Soc.} {\bf 13} (1883) 36.}
  \ref{Cayley}{Cayley,A. {\it Quart.Math.J.} {\bf 7} (1866) 304.}
  \ref{Seade}{Seade,J.A. {\it Anal.Inst.Mat.Univ.Nac.Aut\'on
  M\'exico} {\bf 21} (1981) 129.}
  \ref{CM}{Cisneros--Molina,J.L. {\it Geom.Dedicata} {\bf84} (2001)
  \ref{Goette1}{Goette,S. \jram {526} {2000} 181.}
  207.}
  \ref{NandO}{Nash,C. and O'Connor,D--J, \jmp {36}{1995}{1462}.}
  \ref{Dows}{Dowker,J.S. \aop{71}{1972}{577}; Dowker,J.S. and Pettengill,D.F.
  \jpa{7}{1974}{1527}; J.S.Dowker in {\it Quantum Gravity}, edited by
  S. C. Christensen (Hilger,Bristol,1984)}
  \ref{Jadhav2}{Jadhav,S.P. \prD{43}{1991}{2656}.}
  \ref{Dow11}{Dowker,J.S. \cqg{21}{2004}4247.}
  \ref{Dow12}{Dowker,J.S. \cqg{21}{2004}4977.}
  \ref{Dow13}{Dowker,J.S. \jpa{38}{2005}1049.}
  \ref{Zagier}{Zagier,D. \ma{202}{1973}{149}}
  \ref{RandG}{Rademacher, H. and Grosswald,E. {\it Dedekind Sums},
  (Carus, MAA, 1972).}
  \ref{Berndt7}{Berndt,B, \aim{23}{1977}{285}.}
  \ref{HKMM}{Harvey,J.A., Kutasov,D., Martinec,E.J. and Moore,G.
  {\it Localised Tachyons and RG Flows}, hep-th/0111154.}
  \ref{Beck}{Beck,M., {\it Dedekind Cotangent Sums}, {\it Acta Arithmetica}
  {\bf 109} (2003) 109-139 ; math.NT/0112077.}
  \ref{McInnes}{McInnes,B. {\it APS instability and the topology of the brane
  world}, hep-th/0401035.}
  \ref{BHS}{Brevik,I, Herikstad,R. and Skriudalen,S. {\it Entropy Bound for the
  TM Electromagnetic Field in the Half Einstein Universe}; hep-th/0508123.}
  \ref{BandO}{Brevik,I. and Owe,C.  \prD{55}{4689}{1997}.}
  \ref{Kenn}{Kennedy,G. Thesis. University of Manchester 1978.}
  \ref{KandU}{Kennedy,G. and Unwin S. \jpa{12}{L253}{1980}.}
  \ref{BandO1}{Bayin,S.S.and Ozcan,M.
  \prD{48}{2806}{1993}; \prD{49}{5313}{1994}.}
  \ref{Chang}{Chang, P., {\it Quantum Field Theory on Regular Polytopes}.
   Thesis. University of Manchester, 1993.}
  \ref{Barnesa}{Barnes,E.W. {\it Trans. Camb. Phil. Soc.} {\bf 19} (1903) 374.}
  \ref{Barnesb}{Barnes,E.W. {\it Trans. Camb. Phil. Soc.}
  {\bf 19} (1903) 426.}
  \ref{Stanley1}{Stanley,R.P. \joa {49Hilf}{1977}{134}.}
  \ref{Stanley}{Stanley,R.P. \bams {1}{1979}{475}.}
  \ref{Hurley}{Hurley,A.C. \pcps {47}{1951}{51}.}
  \ref{IandK}{Iwasaki,I. and Katase,K. {\it Proc.Japan Acad. Ser} {\bf A55}
  (1979) 141.}
  \ref{IandT}{Ikeda,A. and Taniguchi,Y. {\it Osaka J. Math.} {\bf 15} (1978)
  515.}
  \ref{GandM}{Gallot,S. and Meyer,D. \jmpa{54}{1975}{259}.}
  \ref{Flatto}{Flatto,L. {\it Enseign. Math.} {\bf 24} (1978) 237.}
  \ref{OandT}{Orlik,P and Terao,H. {\it Arrangements of Hyperplanes},
  Grundlehren der Math. Wiss. {\bf 300}, (Springer--Verlag, 1992).}
  \ref{Shepler}{Shepler,A.V. \joa{220}{1999}{314}.}
  \ref{SandT}{Solomon,L. and Terao,H. \cmh {73}{1998}{237}.}
  \ref{Vass}{Vassilevich, D.V. \plb {348}{1995}39.}
  \ref{Vass2}{Vassilevich, D.V. \jmp {36}{1995}3174.}
  \ref{CandH}{Camporesi,R. and Higuchi,A. {\it J.Geom. and Physics}
  {\bf 15} (1994) 57.}
  \ref{Solomon2}{Solomon,L. \tams{113}{1964}{274}.}
  \ref{Solomon}{Solomon,L. {\it Nagoya Math. J.} {\bf 22} (1963) 57.}
  \ref{Obukhov}{Obukhov,Yu.N. \pl{109B}{1982}{195}.}
  \ref{BGH}{Bernasconi,F., Graf,G.M. and Hasler,D. {\it The heat kernel
  expansion for the electromagnetic field in a cavity}; math-ph/0302035.}
  \ref{Baltes}{Baltes,H.P. \prA {6}{1972}{2252}.}
  \ref{BaandH}{Baltes.H.P and Hilf,E.R. {\it Spectra of Finite Systems}
  (Bibliographisches Institut, Mannheim, 1976).}
  \ref{Ray}{Ray,D.B. \aim{4}{1970}{109}.}
  \ref{Hirzebruch}{Hirzebruch,F. {\it Topological methods in algebraic
  geometry} (Springer-- Verlag,\break  Berlin, 1978). }
  \ref{BBG}{Bla\v{z}i\'c,N., Bokan,N. and Gilkey, P.B. {\it Ind.J.Pure and
  Appl.Math.} {\bf 23} (1992) 103.}
  \ref{WandWi}{Weck,N. and Witsch,K.J. {\it Math.Meth.Appl.Sci.} {\bf 17}
  (1994) 1017.}
  \ref{Norlund}{N\"orlund,N.E. \am{43}{1922}{121}.}
  \ref{Duff}{Duff,G.F.D. \aom{56}{1952}{115}.}
  \ref{DandS}{Duff,G.F.D. and Spencer,D.C. \aom{45}{1951}{128}.}
  \ref{BGM}{Berger, M., Gauduchon, P. and Mazet, E. {\it Lect.Notes.Math.}
  {\bf 194} (1971) 1. }
  \ref{Patodi}{Patodi,V.K. \jdg{5}{1971}{233}.}
  \ref{GandS}{G\"unther,P. and Schimming,R. \jdg{12}{1977}{599}.}
  \ref{MandS}{McKean,H.P. and Singer,I.M. \jdg{1}{1967}{43}.}
  \ref{Conner}{Conner,P.E. {\it Mem.Am.Math.Soc.} {\bf 20} (1956).}
  \ref{Gilkey2}{Gilkey,P.B. \aim {15}{1975}{334}.}
  \ref{MandP}{Moss,I.G. and Poletti,S.J. \plb{333}{1994}{326}.}
  \ref{BKD}{Bordag,M., Kirsten,K. and Dowker,J.S. \cmp{182}{1996}{371}.}
  \ref{RandO}{Rubin,M.A. and Ordonez,C. \jmp{25}{1984}{2888}.}
  \ref{BaandD}{Balian,R. and Duplantier,B. \aop {112}{1978}{165}.}
  \ref{Kennedy2}{Kennedy,G. \aop{138}{1982}{353}.}
  \ref{DandKi2}{Dowker,J.S. and Kirsten, K. {\it Analysis and Appl.}
 {\bf 3} (2005) 45.}
  \ref{Dow40}{Dowker,J.S. \cqg{23}{2006}{1}.}
  \ref{BandHe}{Br\"uning,J. and Heintze,E. {\it Duke Math.J.} {\bf 51} (1984)
   959.}
  \ref{Dowl}{Dowker,J.S. {\it Functional determinants on M\"obius corners};
    Proceedings, `Quantum field theory under
    the influence of external conditions', 111-121,Leipzig 1995.}
  \ref{Dowqg}{Dowker,J.S. in {\it Quantum Gravity}, edited by
  S. C. Christensen (Hilger, Bristol, 1984).}
  \ref{Dowit}{Dowker,J.S. \jpa{11}{1978}{347}.}
  \ref{Kane}{Kane,R. {\it Reflection Groups and Invariant Theory} (Springer,
  New York, 2001).}
  \ref{Sturmfels}{Sturmfels,B. {\it Algorithms in Invariant Theory}
  (Springer, Vienna, 1993).}
  \ref{Bourbaki}{Bourbaki,N. {\it Groupes et Alg\`ebres de Lie}  Chap.III, IV
  (Hermann, Paris, 1968).}
  \ref{SandTy}{Schwarz,A.S. and Tyupkin, Yu.S. \np{242}{1984}{436}.}
  \ref{Reuter}{Reuter,M. \prD{37}{1988}{1456}.}
  \ref{EGH}{Eguchi,T. Gilkey,P.B. and Hanson,A.J. \prp{66}{1980}{213}.}
  \ref{DandCh}{Dowker,J.S. and Chang,Peter, \prD{46}{1992}{3458}.}
  \ref{APS}{Atiyah M., Patodi and Singer,I.\mpcps{77}{1975}{43}.}
  \ref{Donnelly}{Donnelly.H. {\it Indiana U. Math.J.} {\bf 27} (1978) 889.}
  \ref{Katase}{Katase,K. {\it Proc.Jap.Acad.} {\bf 57} (1981) 233.}
  \ref{Gilkey3}{Gilkey,P.B.\invm{76}{1984}{309}.}
  \ref{Degeratu}{Degeratu.A. {\it Eta--Invariants and Molien Series for
  Unimodular Groups}, Thesis MIT, 2001.}
  \ref{Seeley}{Seeley,R. \ijmp {A\bf18}{2003}{2197}.}
  \ref{Seeley2}{Seeley,R. .}
  \ref{melrose}{Melrose}
  \ref{berard}{B\'erard,P.}
  \ref{gromes}{Gromes,D.}
  \ref{Ivrii}{Ivrii}
  \ref{DandW}{Douglas,R.G. and Wojciekowski,K.P. \cmp{142}{1991}{139}.}
  \ref{Dai}{Dai,X. \tams{354}{2001}{107}.}
  \ref{Kuznecov}{Kuznecov}
  \ref{DandG}{Duistermaat and Guillemin.}
  \ref{PTL}{Pham The Lai}
\end{putreferences}

\bye